\UseRawInputEncoding
\documentclass[11pt]{amsart}
\usepackage{mathrsfs}
\usepackage{amssymb}

\usepackage{amscd}
\usepackage{amsmath, amssymb}
\usepackage{amsfonts}
\usepackage[colorlinks,linkcolor=blue,citecolor=blue, pdfstartview=FitH]{hyperref}
\usepackage[all]{xy}

\numberwithin{equation}{section}

\theoremstyle{plain}
\newtheorem{thm}{Theorem}[section]
\newtheorem{theorem}[thm]{Theorem}
\newtheorem{lemma}[thm]{Lemma}
\newtheorem{corollary}[thm]{Corollary}
\newtheorem{proposition}[thm]{Proposition}

\theoremstyle{definition}

\newtheorem{example}[thm]{Example}

\newtheorem{defn-thm}[thm]{Definition-Theorem}

\newcommand{\Image}{{ \textrm{Im}\,}}
\newcommand{\tabincell}[2]{\begin{tabular}{@{}#1@{}}#2\end{tabular}}
\newcommand{\C}{{ \mathbb{C} }}
\newcommand{\R}{{ \mathbb{R} }}
\newcommand{\B}{{ \mathfrak{B} }}
\newcommand{\g}{{ \mathfrak{g} }}
\newcommand{\E}{{ \mathfrak{E} }}
\newcommand{\p}{{ \partial }}
\newcommand{\pb}{{ \bar{\partial} }}
\newcommand{\ot}{{\,\otimes\,}}
\newcommand{\wg}{{\,\wedge\,}}
\newcommand{\Hom}{{ \textrm{Hom}\,}}
\newcommand{\ad}{{ \textrm{ad}\,}}

\begin{document}

\title{Deformations of Dolbeault cohomology classes for Lie algebra with complex structures}
\author{Wei Xia}
\address{Wei Xia, Mathematical Science Research Center, Chongqing University of Technology, Chongqing, P.R.China, 400054.} \email{xiawei@cqut.edu.cn, xiaweiwei3@126.com}

\thanks{This work was supported by the National Natural Science Foundation of China No. 11901590.}
\date{\today}

\begin{abstract}
In this paper, we study deformations of complex structures on Lie algebras and its associated deformations of Dolbeault cohomology classes. A complete deformation of complex structures is constructed in a way similar to the Kuranishi family. The extension isomorphism is shown to be valid in this case. As an application, we prove that given a family of left invariant deformations $\{M_t\}_{t\in B}$ of a compact complex manifold $M=(\Gamma\setminus G, J)$ where $G$ is a Lie group, $\Gamma$ a sublattice and $J$ a left invariant complex structure, the set of all $t\in B$ such that the Dolbeault cohomology on $M_t$ may be computed by left invariant tensor fields is an analytic open subset of $B$.

\vskip10pt
\noindent
{\bf MSC~Classification (2020):} 32G05, 57T15, 53C15, 17B56
\end{abstract}
\maketitle

\section{Introduction}
Let $\mathfrak{g}$ be a real Lie algebra of even dimension, a complex structure on $\g$ is given by a endomorphism $J$ on $\g$ with $J^2=-1$ and such that
\[
[JX,JY]-[X,Y]-J[JX,Y]-J[X,JY]=0,~\forall X,Y\in \g.
\]
In \cite{GT93}, Gigante-Tomassini studied deformations of complex structures on real Lie algebras. They have defined some cohomology groups to show a rigidity result under certain cohomological conditions. In this work, we consider holomorphic deformations of complex structures on real Lie algebras. Compared to \cite{GT93} our approach follows more closely to Kodaira-Spencer's deformation theory of complex analytic structures on compact complex manifolds. In particular, the cohomology groups used in this paper are the standard Lie algebra Dolbeault cohomology~\cite{Rol09,MPPS06,CFP06}.

Let $B$ be an analytic subset in the unit polydisc $\Delta^m\subseteq \C^m$ with $0\in B$ and $(\g,J)$ a Lie algebra with complex structure, by a \emph{holomorphic deformation} of $(\g,J)$ over $B$ we mean a family $\phi (t)\in \g^{0,1*}\otimes\g^{1,0}$ ($t\in B$) such that
\begin{itemize}
  \item[1.] $\phi(t)$ is holomorphic in $t$ and $\phi(0)=0$;
  \item[2.] $\pb\phi (t) =\frac{1}{2}[\phi(t),\phi(t)],~\forall t\in B$,
\end{itemize}
where $\g^{1,0}\subset\g\ot_\R\C$ is the $\sqrt{-1}$-eigenspace of $J$. A (small) holomorphic deformation $\{\phi (t)\}_{t\in B}$ of $(\g,J)$ is said to be \emph{complete} if for any small deformation $\{\psi (s)\}_{s\in D}$ of $(\g,J)$ there exists holomorphic map $h: (D,0)\to (B,0)$ such that $\{\psi (s)\}_{s\in D}$ is equivalent to the pull back of $\{\phi (t)\}_{t\in B}$ by $h$. Our first main result is the following
\begin{theorem}\label{thm-complete-deformation-0}
Let $(\g,J)$ be a Lie algebra with complex structure. Then there is a complete deformation of $(\g,J)$.
\end{theorem}
The construction of this complete deformation is similar to that of the Kuranishi family in classical deformation theory of complex manifolds. A major problem in this case is that Hodge theory is not available in general~\cite{Rol09}. Instead, we have to employ the theory of Moore-Penrose inverses~\cite{Gro77,WWQ18}.

In \cite{Xia19dDol}, we have studied deformations of Dolbeault cohomology classes in the context of deformations of complex manifolds, see \cite{LSY09,LRY15,ZR13,RZ15,RZ18,LZ18,RWZ19} for some related works. Let $X$ be a compact complex manifold and $X_t$ a small deformation of $X$ whose complex structure is represented by $\phi = \phi(t) \in A^{0,1}(X, T^{1,0})$, then we have~\cite[Thm.\,4.4]{Xia19dDol}
\[
H^{p,q}_{\pb_{\phi}}(X) \cong H^{p,q}_{\bar{\partial}_{t}}(X_t)~,~\forall p,~q \geq 0 ,
\]
where $\pb_{\phi}:=\pb-\mathcal{L}_{\phi}^{1,0}=\pb-(i_\phi\p-\p i_\phi)$ and $H^{p,q}_{\pb_{\phi}}(X)=\frac{\ker\pb_{\phi}\cap A^{p,q}(X)}{\Image\pb_{\phi}\cap A^{p,q}(X)}$. We find that similar phenomenon also occur in the present situation. In fact, we have
\begin{theorem}\label{thm-ext-iso-0}
Let $(\g,J_t)$ be a small deformation of $(\g,J)$ where $J_t$ is a complex structure on $\g$ determined by some $\phi=\phi (t)\in \Lambda^{0,1}\otimes\g^{1,0}$. Assume $\E$ is a $\g^{0,1}$-module formed by tensor products or wedge products of $\g^{1,0}$ and $\g^{1,0*}$ and $\E_t$ is the corresponding $\g_t^{0,1}$-module. Then we have
\begin{equation*}
H^{0,\bullet}_{\pb_t}(\E_t)\cong H^{0,\bullet}_{\pb_{\phi}}(\E)~,
\end{equation*}
where $H^{0,\bullet}_{\pb_t}(\E_t)$ is the cohomology of the complex $(\Lambda_{\g_{t}}^{0,\bullet}\ot \E_t,\pb_t)$ and $H^{0,\bullet}_{\pb_{\phi}}(\E)$ the cohomology of $(\Lambda^{0,\bullet}\ot \E,\pb_{\phi})$ with $\Lambda^{p,q}:=\wedge^{p}\mathfrak{g}^{1,0*}\ot \wedge^{q}\mathfrak{g}^{0,1*}$.
\end{theorem}
In view of this result, the deformation theory developed in \cite{Xia19dDol} can be analogously established, see Section \ref{Deformations of Dolbeault cohomology classes}.

The problem of representing Dolbeault cohomology by invariant tensor fields on manifolds of the form $M=(\Gamma\backslash G,J)$, where $G$ is a real Lie group with a lattice $\Gamma\subset G$ and $J$ is a left invariant complex structure on $G$, is extensively studied in recent years, see e.g.~\cite{CFGU00,CF01,CFP06,Con06,MPPS06,Rol09,RTW20} and the references therein. It is believed that a Nomizu type theorem~\cite{Nom54} should hold for the Dolbeault cohomology on complex nilmanifolds~\cite{Rol11b,CFGU00,CF01,FRR19}. As an application of Theorem \ref{thm-ext-iso-0}, we prove the following
\begin{theorem}\label{thm-main-result-0}
Let $M:=(\Gamma\backslash G,J)$ be a compact complex manifold where $G$ is a real Lie group with a lattice $\Gamma\subset G$ and $J$ is a left invariant complex structure on $G$. Denote by $\g$ the Lie algebra of $G$. Assume $E$ is a holomorphic tensor bundle on $M$. Let $\pi: (\mathcal{M}, M)\to (B,0)$ be a small deformation of $M$ such that each fiber $M_t$ of $\pi$ is represented by a $G$-invariant $\phi (t)\in \Lambda^{0,1}\ot\g^{1,0}$ and $E_t$ is the holomorphic tensor bundle on $M_t$ corresponding to $E$, then the set
\[
\left\{t\in B\mid H^{0,\bullet}_{\pb_t}(\E_t)\cong H^{0,\bullet}_{\pb_t}(M_t, E_t)\right\}
\]
is an analytic open subset (i.e. complement of analytic subset) of $B$, where $\E_t$ is the $\g^{0,1}$-module corresponding to $E_t$.
\end{theorem}
If $E=\Omega^p$, it was shown by Console-Fino~\cite[Thm.\,1]{CF01} that the above set is open in Euclidean topology. Thus Theorem \ref{thm-main-result-0} is an improvement of their result. See also \cite{Rol09large,Ang13,Kas13,Kas16,CFK16,AK17b,OV20} for some related works.

Let $\mathcal{C}(\g)$ be the set of all complex structures on $\g$, then in the notations of Theorem \ref{thm-main-result-0} $\mathcal{C}(\g)$ is identical with the set of all left invariant complex structures on $M=(\Gamma\backslash G,J)$. It is well-known that $\mathcal{C}(\g)$ may be viewed as an algebraic subset of the Grassman manifold $Gr(2n;n)$ (see e.g.~\cite{GT93}). An immediate consequence of Theorem \ref{thm-main-result-0} is the following
\begin{corollary}
Under the same conditions as in Theorem \ref{thm-main-result-0}. The set consisting of all left invariant complex structures on $M=(\Gamma\backslash G,J)$ for which the Dolbeault cohomology may be computed by invariant tensor fields is a Zariski open subset of $\mathcal{C}(\g)$.
\end{corollary}

\section{Deformations of complex structures on Lie algebras} \label{Deformations of complex structures on Lie algebras}
\subsection{Lie algebra with complex structures}\label{Lie-alg-with-complex-str}
Let $\mathfrak{g}$ be a real Lie algebra of dimension $2n$, an \emph{almost complex structure} on $\mathfrak{g}$ is by definition a homomorphism $J:\mathfrak{g}\to \mathfrak{g}$ with $J^2=-1$ and such $J$ is called a \emph{complex structure}
\begin{equation}\label{eq-Nijenjuis condition}
[JX,JY]-[X,Y]-J[JX,Y]-J[X,JY]=0,~\forall X,Y\in \mathfrak{g},
\end{equation}
or equivalently
\begin{equation}\label{eq-NN}
[\mathfrak{g}^{1,0},\mathfrak{g}^{1,0}]\subseteq \mathfrak{g}^{1,0},
\end{equation}
where $\mathfrak{g}^{1,0}$ is the $\sqrt{-1}$-eigenspace of $J$ on the complexification $\mathfrak{g}_{\mathbb{C}}=\mathfrak{g}\otimes \mathbb{C}$ and $\mathfrak{g}_{\mathbb{C}}=\mathfrak{g}^{1,0}\oplus\mathfrak{g}^{0,1}$. Such a pair will be called a \emph{Lie algebra with complex structure} and is denoted by $(\mathfrak{g},J)$. Two complex structures $J_0$ and $J_1$ on $\g$ is said to be \emph{equivalent} if there exists a Lie algebra isomorphism $f:\g\to\g$ such that $fJ_0=J_1f$. On the exterior algebra $\bigwedge \mathfrak{g}^{*}=\bigoplus_k\bigwedge^{k} \mathfrak{g}^{*}$, there is the \emph{Chevalley-Eilenberg complex}
\begin{equation*}
\xymatrix@C=0.5cm{
(\bigwedge \mathfrak{g}^{*},d):~0 \ar[r] & \mathfrak{g}^{*} \ar[r]^{d} & \bigwedge^{2} \mathfrak{g}^{*} \ar[r] & \cdots
 \ar[r] & \bigwedge^{2n-1} \mathfrak{g}^{*} \ar[r]^{d} & \bigwedge^{2n} \mathfrak{g}^{*} \ar[r] & 0,}
\end{equation*}
where the differential is given by
\[
d\omega(x,y):= -\omega[x,y],~\forall \omega\in \mathfrak{g}^{*},~x,y\in \mathfrak{g}~.
\]
The complex structure $J$ on $\mathfrak{g}$ give rise to a decomposition
\[
d=\partial + \bar{\partial}~,
\]
such that $\partial^2=\bar{\partial}^2=0$. Given any basis $\{z^\alpha\}_{\alpha=1}^n$ of $\mathfrak{g}^{1,0*}$ and its dual basis $\{z_\alpha\}_{\alpha=1}^n$ for $\mathfrak{g}^{1,0}$ we have
\begin{equation}\label{eq-p-barp}
\partial z^\alpha= -\frac{1}{2}z^\alpha[z_\beta,z_\gamma]z^{\beta\gamma},~\bar{\partial} z^\alpha= -z^\alpha[z_\beta,\bar{z}_\gamma]z^{\beta\bar{\gamma}},~\bar{\partial} z_\alpha= -z^\gamma[\bar{z}_\beta,z_\alpha]\bar{z}^\beta\otimes z_\gamma,
\end{equation}
where $z^{\beta\gamma}:=z^\beta\wedge z^\gamma$ and $z^{\beta\bar{\gamma}}:=z^\beta\wedge z^{\bar{\gamma}}$. Set $\Lambda^{p,q}:=\wedge^{p}\mathfrak{g}^{1,0*}\ot \wedge^{q}\mathfrak{g}^{0,1*}$, the cohomology of the complex $(\Lambda^{p,\bullet} ,\bar{\partial})$ is called the \emph{Dolbeault cohomology} of $(\mathfrak{g},J)$.

\subsection{Deformations of complex structures} \label{Deformations of complex structures}
Given an almost complex structure $J$ on $\g$, any small deformation $J_t$ of $J$ may be represented by an unique $\phi\in \Lambda^{0,1}\otimes\g^{1,0}$. Indeed, let $\{z^\alpha\}_{\alpha=1}^n$ be a basis of $\mathfrak{g}^{1,0*}$ and $\{z_\alpha\}_{\alpha=1}^n$ its dual basis. Likewise, let $\{w^i\}_{i=1}^n$ be a basis of $\mathfrak{g}_{J_t}^{1,0*}$ and $\{w_i\}_{i=1}^n$ its dual basis. Write
\begin{equation}\label{w to z}
\left\{
\begin{array}{ll}
 w^i=A_{\alpha}^i z^\alpha+B_{\alpha}^i \bar{z}^\alpha  \\
 w_i=c^{\alpha}_i z_\alpha+d^{\alpha}_i \bar{z}_\alpha ~, & \\
\end{array} \right.
\end{equation}
then $\phi\in \Lambda^{0,1}\otimes\g^{1,0}$ is defined as
\begin{equation}\label{eq-phi}
\phi:=\phi_\beta^\alpha \bar{z}^\beta\otimes z_\alpha=(A^{-1})_i^{\alpha}B_{\beta}^i \bar{z}^\beta\otimes z_\alpha = -(\bar{c}^{-1})^i_{\beta}\bar{d}^{\alpha}_i \bar{z}^\beta\otimes z_\alpha~.
\end{equation}
Conversely, given a (small) $\phi\in \Lambda^{0,1}\otimes\g^{1,0}$ the deformation $J_t$ of $J$ is determined as follows
\[
\g_{J_t}^{0,1}:=(1-\phi)\g^{0,1},~\text{and}~\g_{J_t}^{1,0}:=(1-\bar{\phi})\g^{1,0}~.
\]
If $J$ is integrable, i.e. it is a complex structure, then it can be checked that
\begin{equation}\label{eq-Maurer-Cartan-equivalence}
[\g_{J_t}^{0,1},\g_{J_t}^{0,1}]\subseteq \g_{J_t}^{0,1} ~~~~\Longleftrightarrow~~~~ \pb\phi =\frac{1}{2}[\phi,\phi]~,
\end{equation}
where the \emph{Fr\"{o}licher-Nijenhuis bracket} $[~,~]:\Lambda^{0,1}\otimes\g^{1,0}\times \Lambda^{0,1}\otimes\g^{1,0}\to \Lambda^{0,2}\otimes\g^{1,0}$ is defined as follows: for any $\varphi,\psi\in\g^*$ and $X,Y\in\g$,
\begin{align*}
[\varphi\otimes X, \psi\otimes Y]&=\varphi\wedge\psi\otimes [X,Y] \\
&-( i_Yd\varphi\wedge\psi\otimes X+ i_Xd\psi\wedge\varphi\otimes Y + i_Y\varphi d\psi\otimes X+ i_X\psi d\varphi\otimes Y )~.
\end{align*}
In fact, set $\xi_\beta^\alpha:=A_{\beta}^ic_i^\alpha$ then by \eqref{eq-phi} $A_{\beta}^id_i^\alpha= -\xi_\beta^\gamma\overline{\phi_\gamma^\alpha}$. It follows from \eqref{w to z} that $[\g_{J_t}^{0,1},\g_{J_t}^{0,1}]\subseteq \g_{J_t}^{0,1}$ if and only if
\begin{align*}
0&=(A^{-1})_k^{\tau}w^k[\bar{w}_i,\bar{w}_j]\bar{A}_{\delta}^i\bar{A}_{\mu}^j\bar{z}^{\delta\mu}\otimes z_\tau\\
&=(A^{-1})_k^{\tau}(A_{\alpha}^k z^\alpha+B_{\alpha}^k \bar{z}^\alpha)[\bar{c}^{\alpha}_i \bar{z}_\alpha+ \bar{d}^{\alpha}_i z_\alpha, \bar{c}^{\alpha}_j \bar{z}_\alpha+ \bar{d}^{\alpha}_j z_\alpha]\bar{A}_{\delta}^i\bar{A}_{\mu}^j\bar{z}^{\delta\mu}\otimes z_\tau\\
&=\left(\bar{d}^{\beta}_i\bar{d}^{\gamma}_jz^\tau[z_\beta,z_\gamma] + (\bar{c}^{\beta}_i\bar{d}^{\gamma}_j- \bar{c}^{\beta}_j\bar{d}^{\gamma}_i)z^\tau[\bar{z}_\beta,z_\gamma]\right)\bar{A}_{\delta}^i\bar{A}_{\mu}^j\bar{z}^{\delta\mu}\otimes z_\tau \\
&+ \left(\phi_\alpha^\tau\bar{c}^{\beta}_i\bar{c}^{\gamma}_j\bar{z}^\alpha[\bar{z}_\beta,\bar{z}_\gamma]
+\phi_\alpha^\tau(\bar{c}^{\beta}_i\bar{d}^{\gamma}_j- \bar{c}^{\beta}_j\bar{d}^{\gamma}_i)\bar{z}^\alpha[\bar{z}_\beta,z_\gamma]
\right)\bar{A}_{\delta}^i\bar{A}_{\mu}^j\bar{z}^{\delta\mu}\otimes z_\tau\\
&=\left(\bar{\xi}_\delta^\alpha\phi_\alpha^\beta \bar{\xi}_\mu^\lambda\phi_\lambda^\gamma z^\tau[z_\beta,z_\gamma] -
2\bar{\xi}_\delta^\beta \bar{\xi}_\mu^\lambda\phi_\lambda^\gamma z^\tau[\bar{z}_\beta,z_\gamma] \right)\bar{z}^{\delta\mu}\otimes z_\tau \\
&+ \left(\phi_\alpha^\tau\bar{\xi}_\delta^\beta \bar{\xi}_\mu^\gamma\bar{z}^\alpha[\bar{z}_\beta,\bar{z}_\gamma]
-2\phi_\alpha^\tau\bar{\xi}_\delta^\beta \bar{\xi}_\mu^\lambda\phi_\lambda^\gamma \bar{z}^\alpha[\bar{z}_\beta,z_\gamma]
\right)\bar{z}^{\delta\mu}\otimes z_\tau
\end{align*}
which is equivalent to (since the matrix $\xi_\beta^\alpha$ is nonsingular)
\begin{align*}
0&=\phi_\alpha^\beta \phi_\lambda^\gamma z^\tau[z_\beta,z_\gamma]\bar{z}^{\alpha\lambda}\otimes z_\tau -
2\phi_\lambda^\gamma z^\tau[\bar{z}_\beta,z_\gamma] \bar{z}^{\beta\lambda}\otimes z_\tau \\
&+ \phi_\alpha^\tau\bar{z}^\alpha[\bar{z}_\beta,\bar{z}_\gamma] \bar{z}^{\beta\gamma}\otimes z_\tau
-2\phi_\alpha^\tau\phi_\lambda^\gamma \bar{z}^\alpha[\bar{z}_\beta,z_\gamma]
\bar{z}^{\beta\lambda}\otimes z_\tau\\
&=2\phi_\beta^\alpha\bar{z}^{\beta}\otimes \pb z_\alpha - 2\phi_\alpha^\tau\pb \bar{z}^\alpha\otimes z_\tau
+ \phi_\alpha^\beta \phi_\lambda^\gamma \bar{z}^{\alpha\lambda}\otimes [z_\beta,z_\gamma] -
2\phi_\alpha^\tau\phi_\lambda^\gamma \bar{z}^\alpha[\bar{z}_\beta,z_\gamma]
\bar{z}^{\beta\lambda}\otimes z_\tau\\
&=-2\pb\phi + [\phi,\phi]~£¬
\end{align*}
where we have used the fact that
\begin{equation}\label{eq-bracket}
[\phi,\phi]=\phi_\alpha^\beta \phi_\lambda^\gamma \bar{z}^{\alpha\lambda}\otimes [z_\beta,z_\gamma] -
2\phi_\alpha^\tau\phi_\lambda^\gamma \bar{z}^\alpha[\bar{z}_\beta,z_\gamma]
\bar{z}^{\beta\lambda}\otimes z_\tau~.
\end{equation}

We list some formulas which is analogous to the geometric setting (see \cite{KMS93,LR11,Xia19deri,Xia19dDol}) as follows:
\begin{proposition}\label{prop-formulas}
Let $K,L\in \g^*\otimes \g$ and $\phi,\phi'\in \Lambda^{0,1}\otimes\g^{1,0}$. Set
\begin{align*}
\mathcal{L}_{K}&=[i_K,d]=i_Kd-di_K,\\
\mathcal{L}_{K}^{1,0}&=[i_K,\p]=i_K\p-\p i_K,\\
\mathcal{L}_{K}^{0,1}&=[i_K,\pb]=i_K\pb-\pb i_K,\\
\end{align*}
then
\begin{itemize}
\item[(1)] $[\mathcal{L}_{K},\mathcal{L}_{L}]:=\mathcal{L}_{K}\mathcal{L}_{L}+\mathcal{L}_{L}\mathcal{L}_{K}=\mathcal{L}_{[K,L]}$~;
\item[(2)] $[\mathcal{L}_{K},i_{L}]:=\mathcal{L}_{K}i_{L}-i_{L}\mathcal{L}_{K}=i_{[K,L]}-\mathcal{L}_{i_{L}K}$~;
\item[(3)] $\mathcal{L}_{\phi}^{0,1}=-i_{\pb \phi}$~;
\item[(4)] $[\mathcal{L}_{\phi}^{1,0},\mathcal{L}_{\phi'}^{1,0}]:=
    \mathcal{L}_{\phi}^{1,0}\mathcal{L}_{\phi'}^{1,0}+\mathcal{L}_{\phi'}^{1,0}\mathcal{L}_{\phi}^{1,0}=\mathcal{L}_{[\phi,\phi']}^{1,0}$ ~;
\item[(5)] $[\pb,\mathcal{L}_{\phi}^{1,0}]:=\pb\mathcal{L}_{\phi}^{1,0}+\mathcal{L}_{\phi}^{1,0}\pb=\mathcal{L}_{\pb\phi}^{1,0}$ ~.
\end{itemize}
\end{proposition}
\begin{proof} We may assume $K=\varphi\ot X$ and $L=\psi\ot Y$ with $\varphi,\psi\in \g^*$ and $X,Y\in \g$.
For any $\omega\in \g^*$, first note that
\begin{equation}\label{eq-[X,Y]}
i_{[X,Y]}\omega=-d\omega(X,Y)=i_Xi_Yd\omega,\quad\text{and}\quad  i_{[X,Y]}d\omega=i_Xdi_Yd\omega-i_Ydi_Xd\omega~.
\end{equation}
Indeed, let $\{e_j\}$ be a basis of $\g$ and $\{\omega^j\in\g^*\}$ be its basis and assume $d\omega^i=\sum_{j,k}\theta_{jk}^i\omega^j\wg \omega^k$, then we have
\[
d\omega^i(e_a,e_b)=(\sum_{j,k}\theta_{jk}^i\omega^j\ot \omega^k-\omega^k\ot \omega^j)(e_a,e_b)=\theta_{ab}^i-\theta_{ba}^i=-i_{e_a}i_{e_b}d\omega^i~.
\]
It follows that
\begin{align*}
i_Z(i_Xdi_Yd\omega-i_Ydi_Xd\omega)=&i_{[Z,X]}i_Yd\omega-i_{[Z,Y]}i_Xd\omega\\
=&i_{[[Z,X],Y]-[[Z,Y],X]}\omega\\
=&i_{[Z,[X,Y]]}\omega\\
=&i_Zi_{[X,Y]}d\omega,
\end{align*}
for any $Z\in \g$.

$(1)$ We have
\begin{align*}
&\mathcal{L}_{K}\mathcal{L}_{L}\omega\\
=&\mathcal{L}_{K}(i_Ld-di_L)\omega\\
=&(i_Kd-di_K)(\psi\wg i_Yd\omega-i_Y\omega d\psi)\\
=&\varphi\wg i_Xd(\psi\wg i_Yd\omega-i_Y\omega d\psi)-d\varphi\wg i_X(\psi\wg i_Yd\omega-i_Y\omega d\psi)+\varphi\wg di_X(\psi\wg i_Yd\omega-i_Y\omega d\psi)\\
=&\varphi\wg(i_Xd\psi\wg i_Yd\omega+ \psi\wg i_Xdi_Yd\omega- i_Y\omega d i_Xd\psi)\\
&-d\varphi\wg (i_X\psi i_Yd\omega- i_Xi_Yd\omega \psi- i_Y\omega i_Xd\psi)\\
\end{align*}
and
\begin{align*}
\mathcal{L}_{L}\mathcal{L}_{K}\omega
=&\psi\wg(i_Yd\varphi\wg i_Xd\omega+ \varphi\wg i_Ydi_Xd\omega- i_X\omega d i_Yd\varphi)\\
&-d\psi\wg (i_Y\varphi i_Xd\omega- i_Yi_Xd\omega \varphi- i_X\omega i_Yd\varphi)
\end{align*}
and it follows from \eqref{eq-[X,Y]} that
\begin{align*}
&\mathcal{L}_{[K,L]}\omega=(i_{[K,L]}d+di_{[K,L]})\omega\\
=&\varphi\wg \psi\wg i_{[X,Y]}d\omega- i_Yd\varphi\wg\psi\wg i_Xd\omega- i_Xd\psi\wg\varphi\wg i_Yd\omega- i_Y\varphi d\psi\wg i_Xd\omega- i_X\psi d\varphi\wg i_Yd\omega \\
&+i_{[X,Y]}\omega (d\varphi\wg\psi-\varphi\wg d\psi)- i_X\omega(di_Y d\varphi\wg\psi-i_Yd\varphi\wg d\psi)- i_Y\omega(di_Xd\psi\wg\varphi-i_Xd\psi\wg d\varphi)\\
=&\varphi\wg \psi\wg (i_Xdi_Yd\omega-i_Ydi_Xd\omega)+\varphi\wg(i_Xd\psi\wg i_Yd\omega-i_Xi_Yd\omega d\psi-i_Y\omega di_Xd\psi)\\
&+\psi\wg(i_Yd\varphi\wg i_Xd\omega+i_Xi_Y d\omega d\varphi-i_X\omega di_Y d\varphi)\\
&+d\varphi\wg(-i_X\psi i_Yd\omega+i_Y\omega i_Xd\psi)+d\psi\wg (-i_Y\varphi i_Xd\omega+i_X\omega i_Yd\varphi)=(\mathcal{L}_{K}\mathcal{L}_{L}+\mathcal{L}_{L}\mathcal{L}_{K})\omega~.\\
\end{align*}

$(2)$ We have
\begin{align*}
&(\mathcal{L}_{K}i_{L}-i_{L}\mathcal{L}_{K})\omega\\
=&i_Y\omega(\varphi\wg i_Xd\psi-i_X\psi d\varphi)-\psi\wg (i_Y\varphi i_Xd\omega-\varphi i_Yi_X d\omega-i_X\omega i_Yd\varphi)~,
\end{align*}
and
\begin{align*}
&(i_{[K,L]}-\mathcal{L}_{i_{L}K})\omega\\
=&i_{[X,Y]}\omega \varphi\wg\psi-i_X\omega i_Yd\varphi\wg\psi-i_Y\omega i_Xd\psi\wg\varphi-i_Y\omega i_X\psi d\varphi
-i_Y\varphi\psi\wg i_Xd\omega~.
\end{align*}

$(3)$ Let $\{z^\alpha\}_{\alpha=1}^n$ be a basis of $\mathfrak{g}^{1,0*}$, $\{z_\alpha\}_{\alpha=1}^n$ its dual basis and $\phi=\phi_\alpha^\beta z^{\bar{\alpha}}\ot z_{\beta}$. It is clear that $\mathcal{L}_{\phi}^{0,1}\bar{z}^{\alpha}=0=i_{\pb \phi}\bar{z}^{\alpha}$. On the other hand,
\begin{align*}
\mathcal{L}_{\phi}^{0,1}z^{\alpha}
=&-2z^{\alpha}[z_{\gamma},\bar{z}_{\beta}]\phi_{\lambda}^{\gamma} z^{\overline{\lambda\beta}}
+\phi_\gamma^\alpha \bar{z}^{\gamma}[\bar{z}_{\beta},\bar{z}_{\lambda}]z^{\overline{\beta\lambda}}\\
=&-i_{\pb \phi}z^{\alpha}~.
\end{align*}

$(4)$ We have
\[
[\mathcal{L}_{\phi},\mathcal{L}_{\phi'}]=\mathcal{L}_{[\phi,\phi']}=\mathcal{L}_{[\phi,\phi']}^{1,0}+\mathcal{L}_{[\phi,\phi']}^{0,1},
\]
on the other hand,
\[
[\mathcal{L}_{\phi},\mathcal{L}_{\phi'}]=[\mathcal{L}_{\phi}^{1,0}-i_{\bar{\partial}\phi}, \mathcal{L}_{\phi'}^{1,0}-i_{\bar{\partial}\phi'}]
=[\mathcal{L}_{\phi}^{1,0},\mathcal{L}_{\phi'}^{1,0}] - [\mathcal{L}_{\phi}^{1,0}, i_{\bar{\partial}\phi'}] - [i_{\bar{\partial}\phi}, \mathcal{L}_{\phi'}^{1,0}] + [i_{\bar{\partial}\phi}, i_{\bar{\partial}\phi'}].
\]
By comparing the bi-degree, $(4)$ follows.

$(5)$ It follow from the Jacobi identity that
\[[\bar{\partial}, \mathcal{L}_{\phi}^{1,0}]=[\bar{\partial}, [i_{\phi}, \partial]]=[[\bar{\partial}, i_{\phi}], \partial]+ [i_{\phi}, [\bar{\partial}, \partial]]=[-\mathcal{L}_{\phi}^{0,1}, \partial]=[i_{\bar{\partial}\phi}, \partial]=\mathcal{L}_{\bar{\partial}\phi}^{1,0}.
\]
\end{proof}

Let $B$ be an analytic subset in the unit polydisc $\Delta^m\subseteq \C^m$ with $0\in B$ and $(\g,J)$ be a Lie algebra with complex structure, by a \emph{holomorphic deformation} of $(\g,J)$ over $B$ we mean a family $\phi (t)\in \Lambda^{0,1}\otimes\g^{1,0}$ ($t\in B$) such that
\begin{itemize}
  \item[1.] $\phi(t)$ is holomorphic in $t$ and $\phi(0)=0$;
  \item[2.] $\pb\phi (t) =\frac{1}{2}[\phi(t),\phi(t)],~\forall t\in B$.
\end{itemize}
For small enough $t$, each $\phi(t)$ determines a complex structure $J_t$ on $\g$. We will call $B$ the \emph{base} and $\g_{t}=(\g,J_t)$ the \emph{fiber}, of the holomorphic deformation. In this case, define $e^{i_{\phi (t)}}:=\sum_{k\geq 0} i_{\phi (t)}^k$, then for any $1\leq p\leq\dim_{\R}\g$ we have
\[
e^{i_{\phi (t)}}:\Lambda^{p,0}\longrightarrow \Lambda^{p,0}_{\g_{t}},
\]
where $\Lambda^{p,0}_{\g_{t}}:=\wg^p\g_{J_t}^{1,0^*}$. Two holomorphic deformations $\{\phi (t)\}_{t\in B}$ and $\{\tilde{\phi}(t)\}_{t\in B}$ of $(\g,J)$ is said to be \emph{equivalent} if there exists Lie algebra isomorphism $f_t:\g\to\g$ such that $J_{\tilde{\phi}(t)}=f_t^{-1}J_{\phi(t)}f_t$ and $f_t$ is holomorphic in $t$ (in the sense that for any $x\in \g_\C$ we have $f_t(x):B\to\g_\C$ is holomorphic).
By a \emph{small deformation} of $(\g,J)$, we mean a germ of holomorphic deformation. The fiber $(\g,J_t)$ of a small deformation will also be called a small deformation. Given a small deformation $\{\phi (t)\}_{t\in B}$ of $(\g,J)$ and a germ of holomorphic map $h: (D,0)\to (B,0)$, where $D$ is an analytic subset in the unit polydisc $\Delta^l\subseteq \C^l$ with $0\in D$, the \emph{pull back of $\{\phi (t)\}_{t\in B}$ by $h$} is also a small deformation $\{\psi (s)\}_{s\in D}$ of $(\g,J)$ defined by
\[
\psi(s):=\phi (h(s))\in\Lambda^{0,1}\otimes\g^{1,0},\quad \forall s\in D~.
\]
A small deformation $\{\phi (t)\}_{t\in B}$ of $(\g,J)$ is said to be \emph{complete} if for any small deformation $\{\psi (s)\}_{s\in D}$ of $(\g,J)$ there exists holomorphic map $h: (D,0)\to (B,0)$ such that $\{\psi (s)\}_{s\in D}$ is equivalent to the pull back of $\{\phi (t)\}_{t\in B}$ by $h$.
\subsection{A complete deformation}\label{A complete deformation}
Recall that the operator $\pb:\Lambda^{p,\bullet}\to \Lambda^{p,\bullet+1}$ can be extended to
\[
\pb:\Lambda^{0,\bullet}\ot \g^{1,0}\to \Lambda^{0,\bullet+1}\ot\g^{1,0},
\]
where given a basis $\{z_\alpha\}$ of $\g^{1,0}$ and $\varphi=\varphi^\alpha\ot z_\alpha\in\Lambda^{0,q}\ot \g^{1,0}$,
\begin{equation}
\pb\varphi:=\pb\varphi^\alpha\ot z_\alpha+(-1)^q \varphi^\alpha\wg \pb z_\alpha~.
\end{equation}
The cohomology of the complex $(\Lambda^{0,\bullet}\ot \g^{1,0},\pb)$ will be denoted by $H^{0,\bullet}_\pb(\g^{1,0})$. Given a holomorphic deformation $(\phi (t),B)$ of $(\g,J)$ (over $B=\Delta\subseteq \C$), then
\[
\pb\phi (t) =\frac{1}{2}[\phi(t),\phi(t)]\Rightarrow \pb\dot{\phi}(0)=\pb\phi_1=0 \Rightarrow \dot{\phi}(0)\in H^{0,1}_\pb(\g^{1,0})~,
\]
where $\phi_1=\frac{\p\phi(t)}{\p t}\mid_{t=0}$.
Now following the classical theory of Kodaira-Spencer-Kuranishi we construct a complete deformation of $(\g,J)$. It is well-known that Hodge theory plays a fundamental role in the classical theory. But it should be noticed that in our case Hodge theory as in the usual form does not apply in general~\cite[Prop.\,1.14]{Rol09} (However, this will not be a problem if $\g$ is nilpotent, see \cite[Rem.\,1.15]{Rol09}). For this reason, we will employ the theory of Moore-Penrose inverse of linear operators instead, see \cite{Gro77} and \cite{WWQ18} for excellent introductions of this topic.

Let $T\in\Hom_{\C}(V_1,V_2)$ be a linear operator between two finite dimensional $\C$-vector spaces $V_1$ and $V_2$. Assume we have chosen Hermitian inner products $h_1, h_2$ on $V_1, V_2$, respectively. The \emph{Moore-Penrose inverse} of $T$, denoted by $T^\dag\in\Hom_{\C}(V_2,V_1)$, is defined as follows:
\begin{equation}
T^\dag y=\left\{
\begin{array}{ll}
(T\mid_{(\ker T)^\bot})^{-1}y,           &\quad \forall y\in \Image T,  \\
0,& \forall y\in (\Image T)^\bot, \\
\end{array} \right.
\end{equation}
where $(\Image T)^\bot$ is the orthogonal complement of $\Image T$ such that $V_2=\Image T\oplus (\Image T)^\bot$. Then it is clear from the definitions that
\[
TT^\dag=P_{\Image T}\quad\text{and}\quad T^\dag T=P_{(\ker T)^\bot}~,
\]
where $P_{\Image T}, P_{(\ker T)^\bot}$ are projection operators onto the subspaces $\Image T, (\ker T)^\bot$, respectively. In fact, the Moore-Penrose inverse operator $T^\dag\in\Hom_{\C}(V_2,V_1)$ of $T$ is uniquely characterized by this property. Notice that we have $\ker T^\dag=(\Image T)^\bot$.

Let $h$ be a Hermitian inner product on $\g_{\C}=\g\ot_{\R}\C$, then the space $\Lambda^{0,\bullet}\ot\g^{1,0}=\oplus_q\Lambda^{0,q}\ot\g^{1,0}$ is equipped with the induced Hermitian inner product from $h$. We denote its associated norm on $\Lambda^{0,\bullet}\ot\g^{1,0}$ by $\|\cdot\|$. We say an element $\varphi\in \Lambda^{0,\bullet}\ot\g^{1,0}$ is \emph{small} if $\|\varphi\|$ is small. This is our convention throughout this paper.

Because the Fr\"{o}licher-Nijenhuis bracket $[-,-]:\Lambda^{0,\bullet}\otimes\g^{1,0}\times \Lambda^{0,\bullet}\otimes\g^{1,0}\to \Lambda^{0,\bullet}\otimes\g^{1,0}$ is bilinear, we have
\begin{equation}\label{eq-FN-estimate}
\|[\varphi,\psi]\|\leq C \|\varphi\|\cdot\|\psi\|,\quad \forall \varphi,\psi\in \Lambda^{0,\bullet}\otimes\g^{1,0}~,
\end{equation}
where $C>0$ is a constant depends only on $(\g, J)$.
Set
\[
\mathcal{H}^{0,q}:=\ker(\pb:\Lambda^{0,q}\ot \g^{1,0}\to \Lambda^{0,q+1}\ot\g^{1,0})\cap\ker(\pb^\dag:\Lambda^{0,q}\ot \g^{1,0}\to \Lambda^{0,q-1}\ot\g^{1,0}) ,~q>0,
\]
and $\mathcal{H}^{0,0}:=\ker\pb\cap\g^{1,0}$, then there is a orthogonal direct sum decomposition
\[
\ker(\pb:\Lambda^{0,q}\ot \g^{1,0}\to \Lambda^{0,q+1}\ot\g^{1,0})=\mathcal{H}^{0,q}\oplus\Image(\pb:\Lambda^{0,q-1}\ot\g^{1,0}\to \Lambda^{0,q}\ot\g^{1,0}),
\]
and therefore the natural map $\mathcal{H}^{0,q}\to H^{0,q}_\pb(\g^{1,0})$ is an isomorphism. As a result, we have
\begin{align*}
&\Lambda^{0,q}\ot\g^{1,0}\\
=&\mathcal{H}^{0,q}\oplus\Image(\pb:\Lambda^{0,q-1}\ot\g^{1,0}\to \Lambda^{0,q}\ot\g^{1,0})\oplus\ker(\pb:\Lambda^{0,q}\ot \g^{1,0}\to \Lambda^{0,q+1}\ot\g^{1,0})^\bot,
\end{align*}
in other words, we have for any $0\leq q\leq n$
\begin{equation}\label{eq-Hodge-decomposition}
1=\mathcal{H}^{0,q}+\pb\pb^\dag+\pb^\dag\pb,\quad\text{on}~\Lambda^{0,q}\ot\g^{1,0},
\end{equation}
where by abuse of notations $\mathcal{H}^{0,q}:\Lambda^{0,q}\ot\g^{1,0}\to \mathcal{H}^{0,q}$ is a projection operator and $\pb^\dag$ is the Moore-Penrose inverse operator of $\pb$. Note that on $\Lambda^{0,0}\ot\g^{1,0}=\g^{1,0}$ we have $P_{(\ker\pb)^\bot}=\pb^\dag\pb$. Let $\eta_1,\cdots, \eta_r$ be a basis of $\mathcal{H}^{0,1}$. Set $\phi (t):=\sum_{k=1}^{\infty} \phi_{k}$, where $\phi_{k}$ is a homogenous polynomial in $t$ of degree $k$ with coefficients in $\Lambda^{0,1}\ot \g^{1,0}$ and
\begin{equation}\label{eq-kuranishi-def}
\left\{
\begin{array}{ll}
\phi_1&=\sum_{\nu=1}^{r}\eta_{\nu}t_{\nu},\\
\phi_k&=\frac{1}{2}\sum_{j=1}^k\bar{\partial}^\dag[\phi_j,\phi_{k-j}],\quad \forall k>1~.
\end{array} \right.
\end{equation}
It follows form \eqref{eq-FN-estimate} that
\begin{equation}\label{eq-kuranishi-estimate}
\|\phi_k\|\leq \frac{C}{2}\|\bar{\partial}^\dag\| \sum_{j=1}^k \|\phi_j\|\cdot \|\phi_{k-j}\|,
\end{equation}
where
\[
\|\bar{\partial}^\dag\|:=\max\left\{ \|\bar{\partial}^\dag\varphi\|  \Big|\varphi\in \Lambda^{0,\bullet}\otimes\g^{1,0}~\text{and}~\|\varphi\|=1 \right\}~.
\]
By \eqref{eq-kuranishi-estimate} and a standard argument \cite[pp.\,162]{MK71}, we may conclude that the power series $\phi (t)=\sum_{k=1}^{\infty} \phi_{k}$ converges for small $t$. On the other hand, it is easy to see that $\phi=\phi (t)$ defined by \eqref{eq-kuranishi-def} is the unique power series solution of
\begin{equation}\label{eq-kuranishi-eq}
\phi=\sum_{\nu=1}^{r}\eta_{\nu}t_{\nu}+\frac{1}{2}\pb^\dag[\phi,\phi],
\end{equation}
with $\phi (0)=0$.
\begin{proposition}\label{prop-unique-power-seris-solution}
$1$. For any $\varphi\in \Lambda^{0,1}\ot \g^{1,0}$, if $\pb\varphi  =\frac{1}{2}[\varphi,\varphi]$ and $\pb^\dag\varphi=0$, then we must have
\[
\varphi = \mathcal{H}^{0,1}\varphi + \frac{1}{2}\pb^\dag[\varphi,\varphi] ~.
\]
$2$. For small $\|t\|$, the equation
\begin{equation}\label{eq-Kuranishi-eq}
\varphi = \sum_{\nu=1}^{r}\eta_{\nu}t_{\nu} + \frac{1}{2}\pb^\dag[\varphi,\varphi] ~,
\end{equation}
has a unique small solution given by $\varphi=\phi(t)$ where $\phi(t)$ is defined by \eqref{eq-kuranishi-def}.
\end{proposition}
\begin{proof}$1.$ If $\pb^\dag\varphi=0$, then $\varphi\in (\Image \pb)^\bot$ and by \eqref{eq-Hodge-decomposition} we have
\[
\varphi=\mathcal{H}^{0,1}\varphi+\pb^\dag\pb\varphi=\mathcal{H}^{0,1}\varphi+\frac{1}{2}\pb^\dag[\varphi,\varphi]~.
\]
$2.$ We already know $\phi (t)$ is a solution, it is left to show the uniqueness. In fact, let $\varphi$ be two solutions of \eqref{eq-Kuranishi-eq} and set $\tau=\varphi-\phi(t)$. Then
\begin{align*}
\tau
=&\frac{1}{2}\pb^\dag([\varphi,\varphi]-[\phi(t),\phi(t)])\\
=&\frac{1}{2}\pb^\dag(2[\tau,\phi(t)]+[\tau,\tau]),
\end{align*}
and therefore,
\[
\|\tau\|\leq  \|\pb^\dag\|\cdot \|\tau\| ( \|\phi(t)\| + \frac{1}{2} \|\tau\|)~.
\]
When $\|t\|$ is small, $\|\phi(t)\|$ is also small. In the mean time for small enough $\|\varphi\|$ we will have $\|\pb^\dag\|\cdot( \|\phi(t)\| + \frac{1}{2} \|\tau\|)<1$ and $\|\tau\|< \|\tau\|$ if $\tau\neq 0$. Hence we must have $\tau= 0$.
\end{proof}

\begin{proposition}\label{prop-holo=harmonic part vanish}
Let $\phi(t)=\sum_{k=1}^{\infty} \phi_{k}\in \Lambda^{0,1}\ot \g^{1,0}$ be defined by \eqref{eq-kuranishi-def}, then for small $t$ we have
\[
\pb\phi(t)  =\frac{1}{2}[\phi(t),\phi(t)]  \Leftrightarrow \mathcal{H}^{0,2} [\phi(t),\phi(t)] =0~.
\]
\end{proposition}
\begin{proof}If $\pb\phi(t)  =\frac{1}{2}[\phi(t),\phi(t)]$, then it follows form the orthogonal decomposition of $\Lambda^{0,2}\ot \g^{1,0}$ that $\mathcal{H}^{0,2} [\phi(t),\phi(t)] =0$.

Conversely, assume $\mathcal{H}^{0,2} [\phi(t),\phi(t)] =0$ and set $\psi(t)=\pb\phi(t)-\frac{1}{2}[\phi(t),\phi(t)]$. It follows from \eqref{eq-kuranishi-eq} and \eqref{eq-Hodge-decomposition} that
\begin{align*}
\psi(t)
&=\frac{1}{2}\pb\pb^\dag[\phi(t),\phi(t)]-\frac{1}{2}[\phi(t),\phi(t)]\\
&=-\frac{1}{2}\pb^\dag\pb[\phi(t),\phi(t)]\\
&=-\pb^\dag[\pb\phi(t),\phi(t)]\\
&=-\pb^\dag[\psi(t)+\frac{1}{2}[\phi(t),\phi(t)],\phi(t)]\\
&=-\pb^\dag[\psi(t),\phi(t)]~.
\end{align*}
Hence,
\[
\|\psi(t)\|\leq C\|\pb^\dag\|\cdot \|\psi(t)\| \cdot \|\phi(t)\|~.
\]
Now for small $|t|$ we will have $C\|\pb^\dag\|\cdot\|\phi(t)\|<1$ and $\|\psi(t)\|< \|\psi(t)\|$ which is a contradiction if $\psi(t)\neq 0$. So we must have $\psi(t)=0$ whenever $|t|$ is small enough.
\end{proof}

For small $\epsilon>0$, we set
\begin{equation}\label{eq-kuranishi-space}
\B:=\left\{t=(t_1,\cdots, t_r)\in \mathbb{C}^r\mid |t|<\epsilon~\text{and}~\mathcal{H}^{0,2}[\phi (t),\phi (t)]=0  \right\}~,
\end{equation}
then by Proposition \ref{prop-holo=harmonic part vanish} the family $\phi (t)\in \Lambda^{0,1}\otimes\g^{1,0}$ ($t\in \B$) defined by \eqref{eq-kuranishi-def} and \eqref{eq-kuranishi-space} is a holomorphic deformation of $(\g, J)$ over $\B$. We will call this family the \emph{Kuranishi deformation} of $(\g, J)$.

Given $X\in \g^{1,0}$, we have a natural map
\[
\ad X:\g_{\C}\to \g_{\C},\quad Y\longrightarrow [X,Y]~,
\]
and the following exponential map is well-defined
\begin{equation}\label{eq-exponential}
e^{\ad X}:\g_{\C}\to \g_{\C},\quad Y\longrightarrow Y+[X,Y]+\frac{1}{2}[X,[X,Y]]+\cdots~.
\end{equation}
\begin{lemma}\label{lem-adjoint-action}For any $\xi,\eta\in \g_{\C}$, we have
\[
e^{\ad X}[\xi,\eta]=[e^{\ad X}\xi,e^{\ad X}\eta]~.
\]
\end{lemma}
\begin{proof}It follows from the Jacobi identity that $\ad X[\xi,\eta]=[\ad X\xi,\eta]+[\xi,\ad X\eta]$. We have
\begin{align*}
[\xi,\eta]&=[\xi,\eta]~,\\
\ad X[\xi,\eta]&=[\ad X\xi,\eta]+[\xi,\ad X\eta]~,\\
(\ad X)^2[\xi,\eta]&=[(\ad X)^2\xi,\eta]+2[\ad X\xi,\ad X\eta]+[\xi,(\ad X)^2\eta]~,\\
&\cdots
\end{align*}
which implies
\[
e^{\ad X}[\xi,\eta]=[e^{\ad X}\xi,\eta]+[e^{\ad X}\xi,\ad X\eta]+[e^{\ad X}\xi,\frac{1}{2}(\ad X)^2\eta]+\cdots=[e^{\ad X}\xi,e^{\ad X}\eta].
\]
\end{proof}
We see from this Lemma that $e^{\ad X}$ is a Lie algebra isomorphism. In particular, if $J_0$ is a complex structure on $\g$, then it can be checked directly by using \eqref{eq-Nijenjuis condition} that $e^{-\ad X}J_0 e^{\ad X}$ is also a complex structure.
\begin{theorem}\label{thm-completeness-canonical-deformation}
Let $(\g,J)$ be a Lie algebra with complex structure and $\psi\in \Lambda^{0,1}\ot \g^{1,0}$ a small deformation of $(\g,J)$.

$1.$ If $\pb^\dag\psi=0$, then $\psi=\phi (t)$ for some unique $t\in \B$ where $\{\phi (t)\}_{t\in \B}$ is the Kuranishi deformation of $(\g,J)$.

$2.$ There exists unique small $X\in (\ker\pb)^\bot\cap\g^{1,0}$ such that the deformation $\tilde{\psi}\in \Lambda^{0,1}\ot \g^{1,0}$ corresponds to the complex structure $e^{-\ad X}J_\psi e^{\ad X}$ satisfies
\[
\pb^\dag\tilde{\psi}=0~.
\]
\noindent Furthermore, if $\psi=\psi (s)$ is holomorphic in $s\in \C^l$, then the corresponding $X=X(s)$ is also holomorphic in $s$.
\end{theorem}
\begin{proof}$1.$ If $\pb^\dag\psi=0$, then by the first part of Proposition \ref{prop-unique-power-seris-solution} we have
\[
\psi=\mathcal{H}^{0,1}\psi+\frac{1}{2}\pb^\dag[\psi,\psi]~.
\]
Since $\psi$ is a small deformation, we must have $\mathcal{H}^{0,1}\psi=\sum_{\nu=1}^{r}\eta_{\nu}t_{\nu}$ for some unique small $t\in\C^r$. Then by the second part of Proposition \ref{prop-unique-power-seris-solution} we have $\psi=\phi (t)$ where $\{\phi (t)\}_{t\in \B}$ is the Kuranishi deformation.

$2.$ Let $\{z^\alpha\}_{\alpha=1}^n, \{w^i\}_{i=1}^n$ be basis of $\mathfrak{g}^{1,0*}, \mathfrak{g}^{1,0*}_{J_\psi}$ respectively, and $\{z_\alpha\}_{\alpha=1}^n, \{w_i\}_{i=1}^n$ its dual basis. Write
\[
 w^i=A_{\alpha}^i z^\alpha+B_{\alpha}^i \bar{z}^\alpha,
\]
then $\psi=\psi_\beta^{\alpha}\bar{z}^\beta\otimes z_\alpha=(A^{-1})_i^{\alpha}B_{\beta}^i \bar{z}^\beta\otimes z_\alpha$. Let $X\in(\ker\pb)^\bot\cap\g^{1,0}$. Note that $\{(e^{\ad X})^*w^i\}_{i=1}^n$ is a basis of $\mathfrak{g}^{1,0*}_{e^{-\ad X}J_\psi e^{\ad X}}$ where $(e^{\ad X})^*:\g_{\C}^*\to \g_{\C}^*$ is induced from $e^{\ad X}$. Set $\tilde{A}_{\alpha}^i=(e^{\ad X})^*w^i(z_\alpha)$ and $\tilde{B}_{\alpha}^i=(e^{\ad X})^*w^i(\bar{z}_\alpha)$, we have
\begin{equation}\label{eq-A}
\tilde{A}_{\alpha}^i=w^i(e^{\ad X}z_\alpha)=A_{\alpha}^i z^\alpha(z_\alpha+[X,z_\alpha]+\frac{1}{2}[X,[X,z_\alpha]]+\cdots) ~,
\end{equation}
and
\begin{equation}\label{eq-B}
\tilde{B}_{\alpha}^i=w^i(e^{\ad X}\bar{z}_\alpha)=(A_{\alpha}^i z^\alpha+B_{\alpha}^i \bar{z}^\alpha)(\bar{z}_\alpha+[X,\bar{z}_\alpha]+\frac{1}{2}[X,[X,\bar{z}_\alpha]]+\cdots) .
\end{equation}
Hence,
\begin{align*}
\tilde{\psi}
&=(\tilde{A}^{-1})_i^{\alpha}\tilde{B}_{\beta}^i \bar{z}^\beta\otimes z_\alpha\\
&=(A^{-1})_i^{\lambda}(\delta^\lambda_\alpha-z^\lambda[X,z_\alpha]+R_1)(B_{\beta}^i+A_{\gamma}^iz^\gamma[X,\bar{z}_\beta]+R_2)\bar{z}^\beta\otimes z_\alpha\\
&=\psi+\pb X-(\psi_\beta^{\lambda}z^\lambda[X,z_\alpha]+z^\lambda[X,z_\alpha]z^\lambda[X,\bar{z}_\beta] )\bar{z}^\beta\otimes z_\alpha+ R(\psi, X),
\end{align*}
where $R_1=A_\lambda^i (\tilde{A}^{-1})_i^{\alpha}-\delta^\lambda_\alpha+z^\lambda[X,z_\alpha]$, $R_2=\tilde{B}_{\beta}^i-B_{\beta}^i-A_{\gamma}^iz^\gamma[X,\bar{z}_\beta]$, and
\[
R(\psi, X)=\left((A^{-1})_i^{\alpha}R_2+R_1B_{\beta}^i+R_1A_{\gamma}^i-z^\lambda[X,z_\alpha]R_2+R_1R_2\right)\bar{z}^\beta\otimes z_\alpha~.
\]
In particular, if $\psi=0$, then $B_{\beta}^i=0$ and we may assume $A_{\alpha}^i=\delta_{\alpha}^i$, in which case, it follows from \eqref{eq-A} and \eqref{eq-B} that both $R_1=R_1(X),~R_2=R_2(X)$ are quadratic in $X$:
\[
R_a(sX)=s^2R_a(X),\quad a=1,2,~\forall s\in\C~.
\]
Hence $R(0, X)$ is quadratic in $X$. Now $\pb^\dag\tilde{\psi}=0$ holds if and only if
\[
\pb^\dag\psi+\pb^\dag\pb X-\pb^\dag\left((\psi_\beta^{\lambda}z^\lambda[X,z_\alpha]+z^\lambda[X,z_\alpha]z^\lambda[X,\bar{z}_\beta] )\bar{z}^\beta\otimes z_\alpha\right)+ \pb^\dag R(\psi, X)=0,
\]
or
\[
X+\pb^\dag\psi-\pb^\dag\left((\psi_\beta^{\lambda}z^\lambda[X,z_\alpha]+z^\lambda[X,z_\alpha]z^\lambda[X,\bar{z}_\beta] )\bar{z}^\beta\otimes z_\alpha\right)+ \pb^\dag R(\psi, X)=0,
\]
where we have used the fact $\pb^\dag\pb X=P_{(\ker\pb)^\bot}X=X$.

Consider the map $\Phi:(\Lambda^{0,1}\ot \g^{1,0})\times \g^{1,0}\to \g^{1,0}$ defined by
\[
 \Phi(\psi, X):=X+\pb^\dag\psi-\pb^\dag\left((\psi_\beta^{\lambda}z^\lambda[X,z_\alpha]+z^\lambda[X,z_\alpha]z^\lambda[X,\bar{z}_\beta] )\bar{z}^\beta\otimes z_\alpha\right)+ \pb^\dag R(\psi, X),
\]
which satisfies $\Phi(0,0)=0$ and is holomorphic near $(0,0)$. Since for any $Y\in\g^{1,0}$, we have
\begin{align*}
\frac{\p\Phi}{\p X}_{(0,0)}(Y)&=\frac{d}{ds}\Phi(0,sY)\Big|_{s=0}\\
&=\frac{d}{ds}\left\{sY-s^2\pb^\dag\left((z^\lambda[Y,z_\alpha]z^\lambda[Y,\bar{z}_\beta] )\bar{z}^\beta\otimes z_\alpha\right)+ s^2\pb^\dag R(0, Y)\right\}\Big|_{s=0}=Y,
\end{align*}
in other words, $\frac{\p\Phi}{\p X}\Big|_{(0,0)}:\g^{1,0}\to\g^{1,0}$ is the identity map, the conclusion then follows from the implicit function theorem.

The last assertion follows since the implicit function theorem implies there is a holomorphic map from a neighbourhood of $0\in \Lambda^{0,1}\ot \g^{1,0}$ to a neighbourhood of $0\in\g^{1,0}$ sending $\psi$ to $X$.
\end{proof}

\begin{corollary}\label{coro-complete-deformation}
Let $(\g,J)$ be a Lie algebra with complex structure, its Kuranishi deformation $\{\phi (t)\}_{t\in \B}$ defined by \eqref{eq-kuranishi-def} and \eqref{eq-kuranishi-space} is a complete deformation.
\end{corollary}
\begin{proof}Let $\{\psi (s)\}_{s\in D}$ be a small deformation of $(\g,J)$ where $D$ is an analytic subset in the unit polydisc $\Delta^l\subseteq \C^l$ with $0\in D$. By Theorem \ref{thm-completeness-canonical-deformation}, we have another small deformation $\{\tilde{\psi}(s)\}_{s\in D}$ which is equivalent to $\{\psi (s)\}_{s\in D}$ such that
\[
\tilde{\psi}(s)=\phi(t),~t\in \B.
\]
The correspondence $s\mapsto t$ defines a function $h:D\to \B$ with $h(s)=t,~\forall s\in D$. It follows that
\[
0=\frac{\p \tilde{\psi}}{\p \bar{s}^j}(0)=\sum_{i=1}^r\frac{\p \phi}{\p t^i}(0)\frac{\p h^i}{\p \bar{s}^j}(0)\Rightarrow\frac{\p h^i}{\p \bar{s}^j}(0)=0,\quad i=1,\cdots,r;~j=1,\cdots,l~.
\]
\end{proof}

\subsection{Applications to deformations of complex manifolds}
Let $M:=(\Gamma\backslash G,J)$ be a compact complex manifold where $G$ is a real Lie group with a lattice $\Gamma\subset G$ and $J$ is a left invariant complex structure on $G$. Denote by $\g$ the Lie algebra of $G$. It is will-known $\g$ may be identified with the space of $G$-left invariant vector fields on $G$, then the left invariant complex structure $J$ on $M$ determines a complex structure, which is still denoted by $J$, on $\g$ and we have the decomposition
\[
\g\ot_\R\C=\g^{1,0}\oplus\g^{0,1}~,
\]
where $\g^{1,0}$ is the space of $G$-left invariant vector fields of type $(1,0)$ on $G$. By Corollary \ref{coro-complete-deformation}, the Kuranishi deformation of $(\g, J)$ is consists of all small left invariant deformations of $M=(\Gamma\backslash G,J)$.
The complex $(\Lambda^{0,\bullet}\ot \g^{1,0},\pb)$ when viewed as the space of $G$-invariant tensor fields on $G$, is naturally a subcomplex of $(A^{0,\bullet}(M)\ot T^{1,0}_M,\pb)$ which is regarded as the space of $\Gamma$-invariant tensor fields on $G$. This inclusion map induce a natural map on cohomology
\[
H^{0,\bullet}_\pb(\g^{1,0})\longrightarrow H^{0,\bullet}_\pb(M, T^{1,0}_M)~.
\]
\begin{theorem}\label{thm-left-inv-kuranishi}\cite{MPPS06,CFP06}
Let $M:=(\Gamma\backslash G,J)$ be a compact complex manifold where $G$ is a real Lie group with a lattice $\Gamma\subset G$ and $J$ is a left invariant complex structure on $G$. Denote by $\g$ the Lie algebra of $G$. Assume the natural map
\[
H^{0,1}_\pb(\g^{1,0})\longrightarrow H^{0,1}_\pb(M, T^{1,0}_M)~,
\]
induced by inclusion is surjective, then any small deformation of $M$ is again a left invariant complex structure.
\end{theorem}
\begin{proof}
Let $\{\phi (t)=\sum_k \phi_k\}_{t\in \mathcal{B}}$ be the Beltrami differentials of the Kuranishi family $\pi: (\mathcal{X}, X_0)\to (\mathcal{B},0)$ of $M$, it is enough to show $\phi (t)$ is $G$-invariant for each $t\in\mathcal{B}$. By the definition of Kuranishi family, $\phi_1 =\sum_{\nu=1}^{m}\eta_{\nu}t_{\nu}$ where $\eta_1,\cdots, \eta_r$ is a basis of the harmonic space $\mathcal{H}^{0,1}_\pb(M, T^{1,0}_M)$.

In fact, choose a $G$-invariant Hermitian metric on $M$, then the operators $\pb, \pb^*, G_\pb$ in Hodge theory preserves $G$-invariant tensor fields and we have $\pb^\dag=\pb^*G_\pb=G_\pb\pb^*$ on $\Lambda^{0,\bullet}\ot \g^{1,0}$. It follows from our assumption that the following map
\[
\mathcal{H}^{0,1}\to H^{0,1}_\pb(\g^{1,0})\to H^{0,1}_\pb(M, T^{1,0}_M)\to \mathcal{H}^{0,1}_\pb(M, T^{1,0}_M)
\]
is surjective. Furthermore, since $\ker\pb^\dag=\ker G_\pb\pb^*=\ker \pb^*$ on $\Lambda^{0,\bullet}\ot \g^{1,0}$, we have
\[
\mathcal{H}^{0,1}=\ker\pb\cap \ker\pb^\dag\cap (\Lambda^{0,1}\ot \g^{1,0})=\ker\pb\cap \ker\pb^*\cap (\Lambda^{0,1}\ot \g^{1,0})
\subset\mathcal{H}^{0,1}_\pb(M, T^{1,0}_M).
\]
We see that $\mathcal{H}^{0,1}= \mathcal{H}^{0,1}_\pb(M, T^{1,0}_M)$ and $\phi_1$ is $G$-invariant. As a result, each term $\phi_k$ in $\phi (t)=\sum_k \phi_k$ is $G$-invariant and in particular we have $\phi (t)\in \Lambda^{0,1}\ot \g^{1,0}$ for each $t\in\mathcal{B}$.
\end{proof}

\begin{corollary}
Let $M:=(\Gamma\backslash G,J)$ be a compact complex manifold where $G$ is a real Lie group with a lattice $\Gamma\subset G$ and $J$ is a left invariant complex structure on $G$. Denote by $\g$ the Lie algebra of $G$ and $\pi: (\mathcal{X}, M)\to (\mathcal{B},0)$ the Kuranishi family of $M$. Assume the natural map
\[
H^{0,1}_\pb(\g^{1,0})\longrightarrow H^{0,1}_\pb(M, T^{1,0}_M)~,
\]
induced by inclusion is surjective, then we have $\mathcal{B}=\B$ (as germs of complex spaces) where $\B$ is the base space of the Kuranishi deformation of $(\g, J)$ defined by \eqref{eq-kuranishi-space}.
\end{corollary}

By a \emph{complex nilmanifold}, we mean a compact complex manifold of the form $\Gamma\backslash G$, where $G$ is a (real) simply connected connected nilpotent Lie group with a left invariant complex structure $J$, and $\Gamma$ is a lattice of $G$ of maximal rank. The complex structure on $\Gamma\backslash G$ is inherit from that of $G$. The left invariant complex structure on $G$ also induce a complex structure, still denoted by $J$, on the Lie algebra $\g$ of $G$.

\begin{corollary}\cite[Thm.\,2.6]{Rol09}
Let $M:=(\Gamma\backslash G,J)$ be a complex nilmanifold of complex dimension $n$. Assume the natural map
\[
H^{1,n-1}_\pb(\g, J)\longrightarrow H^{1,n-1}_\pb(M)~,
\]
induced by inclusion is surjective, then any small deformation of $M$ is again a left invariant complex structure.
\end{corollary}
\begin{proof}In this case, by Serre duality \cite[Coro.\,2.6]{Rol09} we have $H^{1,n-1}_\pb(\g, J)\cong H^{0,1}_\pb(\g^{1,0})$ and $H^{1,n-1}_\pb(M)\cong H^{0,1}_\pb(M, T^{1,0}_M\ot K_M)=H^{0,1}_\pb(M, T^{1,0}_M)$ (it is known that the canonical bundle of any complex nilmanifold is trivial, see \cite{BDV09}).
\end{proof}
\section{Deformations of cohomology classes}
\label{Deformations of cohomology classes}
\subsection{The cohomology of modules over a Lie algebra}
First, we recall some definitions and useful facts, see \cite{Wei94} for more details. Let $\g$ be a real Lie algebra, a (left) \emph{$\g$-module} $M$ is a real vector space equipped with a $\R$-bilinear product
\[
\g\ot M\to M: x\ot m\to xm
\]
such that
\[
[x,y]m=x(ym)-y(xm),\quad \forall~x,y\in\g~\text{and}~m\in M~.
\]
Given two $\g$-modules $M, N$, a \emph{$\g$-module homomorphism} $f:M\to N$ is a $\R$-linear map such tat $f(xm)=xf(m)$, where $x\in\g$ and $m\in M$. Given a $\g$-module $M$, its \emph{invariant submodule} $M^{\g}$ is defined by
\[
M^{\g}=\{m\in M\mid~ xm=0,\quad~\forall x\in \g\}~.
\]
The correspondence $-^{\g}:M\mapsto M^{\g}$ is a left exact functor from the category of $\g$-modules to the category of $\R$-vector spaces. The Lie algebra cohomology of the $\g$-module $M$, denoted by $H^\bullet(\g, M)$ , is defined to be the right derived functor of $-^{\g}$, i.e.
\begin{equation}\label{eq-Lie-alg-coho}
H^\bullet(\g, M):=R^\bullet(-^{\g})(M)~,
\end{equation}
and we call $H^\bullet(\g, M)$ the \emph{cohomology groups of $\g$ with coefficients in $M$}. It is well-known that $H^\bullet(\g, M)$ may be computed by the \emph{Chevalley-Eilenberg complex} $(\Hom_{\R}(\bigwedge^\bullet\g, M), \delta)$, where a $q$-cochain $f\in \Hom_{\R}(\bigwedge^q\g, M)$ is an alternating $\R$-multilinear function $f(x_1,\cdots,x_q)$ of $q$ variables in $\g$, taking values in $M$. The differential $\delta$ of this complex is given by
\begin{align*}
\delta f(x_1,\cdots,x_{q+1})=&\sum_{i=1}^{q+1}(-1)^{i+1}x_i f(x_1,\cdots,\hat{x}_{i},\cdots,x_{q+1})\\
&+\sum_{1\leq i<j\leq q+1}(-1)^{i+j}f([x_i,x_j],x_1,\cdots,\hat{x}_{i},\cdots,\hat{x}_{j},\cdots,x_{q+1})~.
\end{align*}
In other words, we have
\[
H^\bullet(\g, M)\cong \frac{\ker \delta \cap \Hom_{\R}(\bigwedge^\bullet\g, M)}{{\Image \delta} \cap \Hom_{\R}(\bigwedge^\bullet\g, M)}~.
\]
When $M=\R$ is regarded as a trivial $\g$-module, this complex reduce to the complex $(\bigwedge^\bullet\g, d)$ considered in Section \ref{Lie-alg-with-complex-str}.

\subsection{The extension isomorphism}
\label{The extension isomorphism}
Now let $(\g, J)$ be a Lie algebra with complex structure, there is a natural $\g^{0,1}$-module structure on $\Lambda^{p,0}=\bigwedge^p\g^{*1,0}$ as follows: for $p=1$,
\begin{align*}
\g^{0,1}\ot \g^{*1,0}\longrightarrow \g^{*1,0}:~\bar{X}\ot \omega\longmapsto i_{\bar{X}}\pb\omega~,
\end{align*}
where $\bar{X}\in\g^{0,1}$ and $\omega\in\g^{*1,0}$; for $p>1$,
\begin{align*}
\g^{0,1}\ot \Lambda^{p,0}\longrightarrow \Lambda^{p,0}:~\bar{X}\ot z^{\alpha_1}\wg\cdots\wg z^{\alpha_p}\longmapsto
\sum_{j=1}^p z^{\alpha_1}\wg\cdots\wg i_{\bar{X}}\pb z^{\alpha_j}\wg \cdots\wg z^{\alpha_p}~,
\end{align*}
where $\{z^{\alpha}\}$ is a basis of $\g^{*1,0}$. In this case, the Lie algebra cohomology $H^{p,\bullet}_{\pb}(\g,J):=H^\bullet(\g^{0,1}, \Lambda^{p,0})$ is called the \emph{Lie algebra Dolbeault cohomology} of $(\g, J)$~\cite{Rol09} and the Chevalley-Eilenberg complex becomes
\begin{equation*}
\xymatrix@C=0.5cm{
(\Lambda^{p,\bullet},\pb):~0 \ar[r] & \Lambda^{p,0} \ar[r]^{\pb} & \Lambda^{p,1} \ar[r] & \cdots
 \ar[r] & \Lambda^{p,n-1} \ar[r]^{\pb} & \Lambda^{p,n} \ar[r] & 0~.}
\end{equation*} Let $(\g,J_t)$ be a small deformation of $(\g,J)$ where $J_t$ is a complex structure on $\g$ determined by some $\phi=\phi (t)\in \Lambda^{0,1}\otimes\g^{1,0}$. Set
\[
\pb_{\phi} = \bar{\partial} - \mathcal{L}_{\phi}^{1,0}:~\Lambda^{p,\bullet}\longrightarrow \Lambda^{p,\bullet+1},
\]
then since by Proposition \ref{prop-formulas} $(\bar{\partial} - \mathcal{L}_{\phi}^{1,0})^2=\mathcal{L}_{-2\pb\phi+[\phi,\phi]}^{1,0}=0$, we have the \emph{deformed Chevalley-Eilenberg complex}
\begin{equation}\label{eq-deformed-Chevalley-E}
\xymatrix@C=0.5cm{
(\Lambda^{p,\bullet},\pb_{\phi}):~0 \ar[r] & \Lambda^{p,0} \ar[r]^{\pb_{\phi}} & \Lambda^{p,1} \ar[r] & \cdots
 \ar[r] & \Lambda^{p,n-1} \ar[r]^{\pb_{\phi}} & \Lambda^{p,n} \ar[r] & 0~,}
\end{equation}
whose cohomology
\[
H^{p,q}_{\pb_{\phi}}(\g, J):=
\frac{\ker \pb_{\phi} \cap \wedge^{p,q}}{{\Image \pb}_{\phi} \cap \wedge^{p,q}}~,
\]
is called the \emph{deformed Dolbeault cohomology}.
\begin{proposition}\label{prop-iso-as-mod}
Let $(\g,J_t)$ be a small deformation of $(\g,J)$ where $J_t$ is a complex structure on $\g$ determined by some $\phi=\phi (t)\in \Lambda^{0,1}\otimes\g^{1,0}$,
then we have the following commutative diagram ($\g^{0,1}_t=\g^{0,1}_{J_t}=(1-\phi)\g^{0,1}$)
\[
\xymatrix{
 \g^{0,1} \ar[d]_{} \bigotimes \Lambda^{p,0} \ar[d]_{(1-\phi)\ot e^{i_\phi}} \ar[r] & \Lambda^{p,0}:&\bar{X}\ot \omega\longmapsto i_{\bar{X}}\pb_{\phi}\omega  \\
 \g^{0,1}_t  \bigotimes \Lambda^{p,0}_{\g_{t}}  \ar[r] & \Lambda^{p,0}_{\g_{t}}\ar[u]^{P_{J}}:&\bar{X}\ot \omega\longmapsto i_{\bar{X}}\pb_{t}\omega~,   }
\]
where $\pb_{t}+\p_{t}=d$ is the decomposition (of $d$) w.r.t. $J_t$ and
\[
P_{J}:\bigwedge\g_{\C}^*=\bigoplus_{i,j}\Lambda^{i,j}\to \Lambda^{p,0}
\]
is the projection operator. Furthermore, the map
\[
\g^{0,1}\otimes \Lambda^{p,0}\to\Lambda^{p,0}:\bar{X}\ot \omega\mapsto i_{\bar{X}}\pb_{\phi}\omega
\]
defines a $\g^{0,1}$-module structure on $\Lambda^{p,0}$ (denoted by $\Lambda^{p,0}_{\phi}$) and
\[
H^\bullet(\g^{0,1}, \Lambda^{p,0}_{\phi})\cong H^{p,q}_{\pb_{\phi}}(\g, J)~.
\]
\end{proposition}
\begin{proof}We only prove it for $p=1$ since the general case is similar. Let $\{z^\alpha\}_{\alpha=1}^n$ be a basis of $\mathfrak{g}^{1,0*}$ and $\{z_\alpha\}_{\alpha=1}^n$ its dual basis. Assume $\phi=\phi^\alpha_\beta \bar{z}^\beta\ot z_\alpha$, we compute
\begin{align*}
&\left(i_{(1-\phi)\bar{z}_\alpha}\pb_{t}(e^{i_\phi}z^\beta)\right) (z_\gamma)\\
=&i_{z_\gamma}i_{(1-\phi)\bar{z}_\alpha}\pb_{t}(e^{i_\phi}z^\beta)\\
=&i_{z_\gamma}i_{(1-\phi)\bar{z}_\alpha}d(e^{i_\phi}z^\beta)\\
=&(e^{i_\phi}z^\beta)[z_\gamma,(1-\phi)\bar{z}_\alpha]\\
=&(1+i_\phi)z^\beta[z_\gamma,\bar{z}_\alpha-\phi_\alpha^\lambda z_\lambda]\\
=&(z^\beta+\phi^\beta_\nu \bar{z}^\nu)[z_\gamma,\bar{z}_\alpha-\phi_\alpha^\lambda z_\lambda]\\
=&z^\beta[z_\gamma,\bar{z}_\alpha]-\phi_\alpha^\lambda z^\beta[z_\gamma,z_\lambda]+\phi^\beta_\nu\bar{z}^\nu[z_\gamma,\bar{z}_\alpha]\\
=&(i_{\bar{z}_\alpha}\pb_{\phi}z^\beta)(z_\gamma)~,
\end{align*}
which implies $P_{J}i_{(1-\phi)\bar{z}_\alpha}\pb_{t}(e^{i_\phi}z^\beta)=i_{\bar{z}_\alpha}\pb_{\phi}(z^\beta)$ for $1\leq\alpha,\beta\leq n$. The last assertion follows since the Chevalley-Eilenberg complex of the $\g^{0,1}$-module $\Lambda^{p,0}_{\phi}$ is just the deformed complex \eqref{eq-deformed-Chevalley-E}.
\end{proof}

\begin{theorem}
Let $(\g,J_t)$ be a small deformation of $(\g,J)$ where $J_t$ is a complex structure on $\g$ determined by some $\phi=\phi (t)\in \Lambda^{0,1}\otimes\g^{1,0}$. Then we have
\begin{equation*}
H^{p,\bullet}_{\pb_t}(\g,J_t)\cong H^{p,\bullet}_{\pb_{\phi}}(\g, J),\quad 0\leq p\leq n~.
\end{equation*}
\end{theorem}
\begin{proof}Since $\g^{0,1}_t=(1-\phi)\g^{0,1}$ and $\Lambda^{p,0}_{\g_{t}}=e^{i_\phi}\Lambda^{p,0}$, Proposition \ref{prop-iso-as-mod} says exactly that the module structures on $\Lambda^{p,0}_{\g_{t}}$ (defined by $\pb_{t}$) and that of $\Lambda^{p,0}$ (defined by $\pb_{\phi}$) are identical. Hence
\[
H^\bullet(\g^{0,1}_t, \Lambda^{p,0}_{\g_{t}})\cong H^\bullet(\g^{0,1}, \Lambda^{p,0}_{\phi})~.
\]
The conclusion follows because $H^{p,\bullet}_{\pb_t}(\g,J_t)=H^\bullet(\g^{0,1}_t, \Lambda^{p,0}_{\g_{t}})$ and $H^\bullet(\g^{0,1}, \Lambda^{p,0}_{\phi})\cong H^{p,\bullet}_{\pb_{\phi}}(\g, J)$.
\end{proof}

The $\g^{0,1}$-module structure on $\g^{1,0}$ is defined similarly:
\[
\g^{0,1} \otimes \g^{1,0} \longrightarrow \g^{1,0}:~\bar{X}\ot Y\longmapsto i_{\bar{X}}\pb Y,
\]
where $\bar{X}\in\g^{0,1}$ and $Y\in\g^{1,0}$. The Chevalley-Eilenberg complex for $\g^{1,0}$ is $(\Lambda^{0,\bullet}\ot \g^{1,0},\pb)$ and we thus have $H^{0,\bullet}_\pb(\g^{1,0})\cong H^{\bullet}(\g^{0,1},\g^{1,0})$. Correspondingly, the deformed Chevalley-Eilenberg complex for $\g^{1,0}$ is $(\Lambda^{0,\bullet}\ot \g^{1,0},~\pb_{\phi})$ where
\[
\pb_{\phi}:\Lambda^{0,\bullet}\ot \g^{1,0}\longrightarrow \Lambda^{0,\bullet+1}\ot \g^{1,0} ~: X\longmapsto \pb X-[\phi,X] ~.
\]
In general, the tensor product or wedge product of $\g^{0,1}$-modules has natural induced $\g^{0,1}$-module structure. If $\E$ is a $\g^{0,1}$-module formed by tensor products or wedge products of $\g^{1,0}$ and $\g^{1,0*}$, we have
\[
H^{0,\bullet}_\pb(\E)\cong H^{\bullet}(\g^{0,1},\E)
\]
where $H^{0,\bullet}_\pb(\E)$ is the cohomology of the complex $(\Lambda^{0,\bullet}\ot \E,\pb)$. The operator $\pb_{\phi}$ can be naturally extended to
\[
\pb_{\phi}:\Lambda^{0,\bullet}\ot \E\longrightarrow \Lambda^{0,\bullet+1}\ot \E
\]
such that it is in compliance with the Leibniz's rule. In fact, we have $\pb_{\phi}=\pb-\langle\phi|$ where
\[
\langle\phi|:\Lambda^{0,\bullet}\ot \E\longrightarrow \Lambda^{0,\bullet+1}\ot \E
\]
is defined by mixing $\mathcal{L}_{\phi}^{1,0}$ and $[\phi,-]$. For example, for $\E=\Lambda^{p,0}\ot \g^{1,0}$ and $\sigma=\varphi\ot\psi\ot X\in\Lambda^{0,q}\ot\Lambda^{p,0}\ot\g^{1,0}$, we have
\[
\langle\phi|\sigma=(\mathcal{L}_{\phi}^{1,0}\varphi)\ot\psi\ot X+(-1)^q\varphi\ot(\mathcal{L}_{\phi}^{1,0}\psi)\ot X+(-1)^{p+q} \varphi\ot\psi\ot [\phi,X],
\]
where $\varphi\in\Lambda^{0,q},~\psi\in\Lambda^{p,0}$ and $X\in\g^{1,0}$. We will often write $\langle\phi|\sigma\rangle:=\langle\phi|\sigma$.
A modification of the above proofs will give rise to following
\begin{theorem}\label{thm-ext-iso}
Let $(\g,J_t)$ be a small deformation of $(\g,J)$ where $J_t$ is a complex structure on $\g$ determined by some $\phi=\phi (t)\in \Lambda^{0,1}\otimes\g^{1,0}$. Assume $\E$ is a $\g^{0,1}$-module formed by tensor products or wedge products of $\g^{1,0}$ and $\g^{1,0*}$ and $\E_t$ is the corresponding $\g_t^{0,1}$-module. Then we have
\begin{equation*}
H^{0,\bullet}_{\pb_t}(\E_t)\cong H^{0,\bullet}_{\pb_{\phi}}(\E)~,
\end{equation*}
where $H^{0,\bullet}_{\pb_t}(\E_t)$ is the cohomology of the complex $(\Lambda_{\g_{t}}^{0,\bullet}\ot \E_t,\pb_t)$ and $H^{0,\bullet}_{\pb_{\phi}}(\E)$ the cohomology of $(\Lambda^{0,\bullet}\ot \E,\pb_{\phi})$.
\end{theorem}
\subsection{Deformations of Dolbeault cohomology classes}\label{Deformations of Dolbeault cohomology classes}
In view of Section \ref{A complete deformation}, much of the theory developed in \cite{Xia19dDol} has an analogue in the present situation which we will sketch.

Let $(\g, J)$ be a Lie algebra with complex structure and $\{\phi (t)\in \Lambda^{0,1}\otimes\g^{1,0}\}_{t\in B}$ a small deformation. Denote by $\E$ a $\g^{0,1}$-module formed by tensor products or wedge products of $\g^{1,0}$ and $\g^{1,0*}$. We will always assume $\g_{\C}=\g\ot_\R\C$ and $\E$ are equipped with Hermitian metrics. Given $[y]\in H_{\bar{\partial}}^{0,q}(\E)$ and $T\subseteq B$, which is a complex subspace of $B$ containing $0$, a \emph{deformation} of $[y]$ (w.r.t. $\{\phi (t)\}_{t\in B}$) on $T$ is a family $\{\sigma (t)\in \Lambda^{0,q}\ot \E\}_{t\in T}$ such that
\begin{itemize}
  \item[1.] $\sigma (t)$ is holomorphic in $t$ and $[\sigma (0)] = [y]\in H_{\bar{\partial}}^{0,q}(\E)$;
  \item[2.] $\bar{\partial}_{\phi(t)}\sigma (t)=\pb\sigma(t)-\langle\phi(t)|\sigma(t)\rangle =0,~\forall t\in T$.
\end{itemize}
If $B$ is smooth and $T=B$ holds, we say $[y]$ has \emph{unobstructed deformation w.r.t. $\{\phi (t)\}_{t\in B}$}. If $T=B$ holds for any $\{\phi (t)\}_{t\in B}$ with smooth $B$, then we say $[y]$ has \emph{unobstructed deformation}. A deformation $\sigma (t)$ of $[y]$ on $T$ is called \emph{canonical} it satisfies
$\sigma(t) = y + \bar{\partial}^\dag \langle\phi(t) | \sigma(t)\rangle$ for any $t\in T$, where $\bar{\partial}^\dag$ is the Moore-Penrose inverse of $\bar{\partial}$. Two deformations $\sigma (t)$ and $\sigma' (t)$ of $[y]$ on $T$ are \emph{equivalent} if
$[\sigma (t) - \sigma' (t)] = 0 \in H^{0,q}_{\bar{\partial}_{\phi(t)}}(\E)$, $\forall t\in T$.

\begin{proposition}\label{prop-12345}Let $\phi=\phi (t)=\sum_k \phi_k$ be a small deformation of $(\g, J)$.

$1$. $\forall \sigma \in \Lambda^{0,q}\ot \E$, if $\bar{\partial}_{\phi} \sigma = \bar{\partial}\sigma - \langle\phi | \sigma \rangle =0$ and $\bar{\partial}^\dag\sigma=0$, then we must have
\[
\sigma = \mathcal{H}^{0,q}\sigma + \bar{\partial}^\dag\langle\phi | \sigma \rangle .
\]

$2$. The equation
\begin{equation}\label{Kuranishi eq}
\sigma=\sigma_0 + \bar{\partial}^\dag \langle\phi(t) | \sigma \rangle ,~~~~~~~~\text{with}~~\sigma_0\in \ker\pb\cap\Lambda^{0,q}\ot \E,
\end{equation}
has a unique small solution given by $\sigma=\sigma(t)=\sum_{k} \sigma_k$, where $\sigma_{k}=\bar{\partial}^\dag \sum_{j=1}^k\langle\phi_j | \sigma_{k-j} \rangle\in \Lambda^{0,q}\ot \E$ for any $k>0$.

$3$. Let $\sigma$ be a solution of the equation \eqref{Kuranishi eq}. Then we have
\[
\bar{\partial}\sigma = \langle\phi(t) | \sigma \rangle  \Leftrightarrow \mathcal{H}^{0,q} \langle\phi(t) | \sigma \rangle =0.
\]

$4$. For any fixed $t$, the following homomorphism
\[
\ker\bar{\partial}\cap (\Lambda^{0,q}\ot \E)\longrightarrow \ker\bar{\partial}^\dag\cap\Image\bar{\partial}_{\phi(t)}\cap (\Lambda^{0,q+1}\ot \E):x_0\longmapsto \bar{\partial}_{\phi(t)}x(t),
\]
is surjective, where $x(t)$ is the unique solution of $x(t)=x_0 + \bar{\partial}^\dag \langle\phi(t) | x(t)\rangle$. Furthermore, we have
\[
\dim H^{0,q}_{\bar{\partial}}(\E)=\dim\left(\ker\bar{\partial}^\dag\cap\ker\bar{\partial}_{\phi(t)}\right)^q+
\dim \left(\ker\bar{\partial}^\dag\cap\Image\bar{\partial}_{\phi(t)}\right)^{q+1},
\]
where $\left(\ker\bar{\partial}^\dag\cap\ker\bar{\partial}_{\phi(t)}\right)^q:=\ker\bar{\partial}^\dag\cap\ker\bar{\partial}_{\phi(t)}\cap (\Lambda^{0,q}\ot \E)$ and $\left(\ker\bar{\partial}^\dag\cap\Image\bar{\partial}_{\phi(t)}\right)^{q+1}:=\ker\bar{\partial}^\dag\cap\Image\bar{\partial}_{\phi(t)}\cap (\Lambda^{0,q+1}\ot \E)$.

$5$. For any fixed $t$, the natural map
\begin{equation}\label{eq-natural map}
\frac{\ker\bar{\partial}^\dag\cap\ker\bar{\partial}_{\phi(t)} \cap (\Lambda^{0,q}\ot \E)}{\ker\bar{\partial}^\dag\cap\Image\bar{\partial}_{\phi(t)} \cap (\Lambda^{0,q}\ot \E)}\longrightarrow H_{\bar{\partial}_{\phi(t)}}^{0,q}(\E)
\end{equation}
is an isomorphism.

$6$. Let $V=\mathbb{C}\{ \sigma_0^1, \cdots, \sigma_0^N \}\subseteq \mathcal{H}^{0,q}(\E)$, we set
\begin{align*}
V_{t}^q:=
&\{ \sum_{l=1}^N a_l\sigma_0^l\in V \mid  (a_1, \cdots, a_N)\in \mathbb{C}^{N}~\text{s.t.}~\bar{\partial}_{\phi(t)}\sigma(t)=0,\\
&\text{where}~\sigma (t)=\sum_{k} \sigma_k~\text{with}~ \sigma_0=\sum_l a_l\sigma_0^l~\text{and}~ \sigma_{k}=\bar{\partial}^\dag\sum_{i+j=k} \langle\phi_i| \sigma_{j} \rangle,~\forall k\neq 0 \},
\end{align*}
and
\begin{align*}
\tilde{f}_t: &V_{t}^q \longrightarrow \ker\bar{\partial}^\dag\cap\ker\bar{\partial}_{\phi(t)}\cap (\Lambda^{0,q}\ot \E),\\
&\sigma_0\longmapsto \sigma(t)=\sum_{k} \sigma_k,~\text{where}~\sigma_{k}=\bar{\partial}^\dag G\sum_{i+j=k} \langle\phi_i | \sigma_{j} \rangle,~\forall k\neq 0,\\
f_t: &V_{t}^q \longrightarrow \frac{\ker\bar{\partial}^\dag\cap\ker\bar{\partial}_{\phi(t)}}{\ker\bar{\partial}^\dag\cap \Image\bar{\partial}_{\phi(t)}}\cong H^{0,q}_{\bar{\partial}_{\phi(t)}}(\E),\quad\sigma_0\longmapsto [\tilde{f}_t(\sigma_0)]~.
\end{align*}
If $V= \mathcal{H}^{0,q}(\E)$, then $\tilde{f}_t$ is an isomorphism and $f_t$ is surjective.
\end{proposition}
\begin{proof}All these statements follows from the same arguments as in \cite{Xia19dDol}.
\end{proof}
An immediate consequence of this proposition are the following results:
\begin{theorem}\label{thm-g-analytic-subset}
Let $(\g, J)$ be a Lie algebra with complex structure and $\{\phi (t)\in \Lambda^{0,1}\otimes\g^{1,0}\}_{t\in B}$ a small deformation of $(\g, J)$. Denote by $\E$ a $\g^{0,1}$-module formed by tensor products or wedge products of $\g^{1,0}$ and $\g^{1,0*}$ and $\E_t$ the corresponding $\g^{0,1}_t$-module where $\g^{0,1}_t:=\left(1-\phi(t)\right)\g^{0,1}$. Then for any nonnegative integer $k$, the set
\[
\left\{t\in B\mid \dim H^{0,\bullet}_{\pb_t}(\E_t)\geq k\right\}
\]
is an analytic subset of $B$.
\end{theorem}
\begin{proof} We may assume $\{\phi (t)\}$ is just the canonical deformation of $(\g, J)$. By Theorem \ref{thm-ext-iso} and Proposition \ref{prop-12345}, we have
\begin{align*}
&\{t\in B\mid \dim H^{0,q}_{\pb_t}(\E_t)\geq k\}\\
=&\{t\in B\mid \dim H^{0,q}_{\pb_{\phi (t)}}(\E)\geq k\}\\
=&\{t\in B\mid \dim V_t^q-\dim\left(\ker\bar{\partial}^\dag\cap\Image\bar{\partial}_{\phi(t)}\right)\geq k\}\\
=&\{t\in B\mid \dim V_t^q+\dim V_t^{q-1}\geq k+\dim H^{0,q-1}_{\bar{\partial}}(\E)\}\\
\end{align*}
Now, let $\{\sigma_0^l\}_{l=1}^N$ be a basis of $V=\mathcal{H}^{0,q}(\E)$, then it follows from the third assertion of Proposition \ref{prop-12345} that
\begin{align*}
V_{t}^q&=\{\sum_l a_l\sigma_0^l\in\mathcal{H}^{0,q}(\E)\mid \sum_{l=1}^N a_l\bar{\partial}_{\phi(t)}\sigma^l(t) =0\}\\
&=\{\sum_l a_l\sigma_0^l\in\mathcal{H}^{0,q}(\E)\mid \sum_{l=1}^N a_l \mathcal{H}^{0,q+1}\langle\phi(t) | \sigma^l(t) \rangle =0\}~,
\end{align*}
where $\sigma^l(t)=\sum_{k} \sigma_k^l$ with $\sigma_{k}^l=\bar{\partial}^\dag\sum_{i+j=k} \langle\phi_i| \sigma_{j}^l \rangle$ for $k> 0$.
Hence we see that
\begin{align*}
&\{t\in B\mid \dim H^{0,q}_{\pb_t}(\E_t)\geq k\}\\
=&\{t\in B\mid \dim V_t^q+\dim V_t^{q-1}\geq k+\dim H^{0,q-1}_{\bar{\partial}}(\E)\}\\
=&\bigcup_{j=0}^N\Big(\{t\in B\mid \dim V_t^q\geq k+\dim H^{0,q-1}_{\bar{\partial}}(\E)-j\}\cap\{t\in B\mid \dim V_t^{q-1}= j\}\Big)\\
=&\bigcup_{j=0}^N\Big(\{t\in B\mid \dim V_t^q\geq k+\dim H^{0,q-1}_{\bar{\partial}}(\E)-j\}\cap\{t\in B\mid \dim V_t^{q-1}\geq j\}\Big)
\end{align*}
is an analytic subset of $B$, where $N=\dim H^{0,q-1}_{\bar{\partial}}(\E)$.
\end{proof}

\begin{theorem}\label{Deformation of Dolbeault cohomology classes: general case}
Let $(\g, J)$ be a Lie algebra with complex structure and $\{\phi (t)\in \Lambda^{0,1}\otimes\g^{1,0}\}_{t\in \B}$ the canonical deformation of $(\g, J)$. Denote by $\E$ a $\g^{0,1}$-module formed by tensor products or wedge products of $\g^{1,0}$ and $\g^{1,0*}$. Let $V=\mathbb{C}\{ \sigma_0^1, \cdots, \sigma_0^N \}$ be a linear subspace of $\mathcal{H}^{0,q}(\E)$ and $\sigma^l(t)=\tilde{f}_t\sigma_{0}^l,~l=1, \cdots, N$. Define an analytic subset $\B(V)$ of $\B$ by
\[
\B(V):=\{t\in \B\mid \mathcal{H}^{0,q} \langle\phi(t) | \sigma^l(t) \rangle = 0, l=1, \cdots, N\},
\]
Then we have
\begin{equation}\label{B(V)-discription}
\B(V)=\{t\in \B\mid \dim V = \dim \Image f_t +\dim \ker f_t\}.
\end{equation}
In particular, we have
\begin{equation}\label{B'-discription}
\B':=\B(\mathcal{H}^{0,q}(\E))=\{t\in \B\mid \dim H^{0,q}_{\pb}(\E)=\dim H^{0,q}_{\pb_t}(\E_t)+\dim \ker f_t\} .
\end{equation}
\end{theorem}
\begin{proof}This follows from $3., 6.$ of Proposition \ref{prop-12345} and Theorem \ref{thm-ext-iso}.
\end{proof}
\begin{theorem}\label{thm-2nd-main}
Let $(\g, J)$ be a Lie algebra with complex structure and $\{\phi (t)\in \Lambda^{0,1}\otimes\g^{1,0}\}_{t\in B}$ a small deformation of $(\g, J)$. Denote by $\E$ a $\g^{0,1}$-module formed by tensor products or wedge products of $\g^{1,0}$ and $\g^{1,0*}$ and $\E_t$ the corresponding $\g^{0,1}_t$-module where $\g^{0,1}_t:=\left(1-\phi(t)\right)\g^{0,1}$. For each $q\geq 0$, set
\[
v^q_t:=\dim H_{\bar{\partial}}^{0,q}(\E)-\dim \ker\bar{\partial}_{\phi(t)}\cap\ker\bar{\partial}^\dag\cap (\Lambda^{0,q}\ot \E),
\]
then we have
\begin{equation}\label{eq-1461}
\dim H_{\bar{\partial}}^{0,q}(\E)=\dim H_{\bar{\partial}_t}^{0,q}(\E_t)+v^q_t+v^{q-1}_t.
\end{equation}
In particular, $\dim H_{\bar{\partial}_t}^{0,q}(\E_t)$ is independent of $t\in B$ if and only if the deformations of classes in $H_{\bar{\partial}}^{0,q}(\E)$ and $H_{\bar{\partial}}^{0,q-1}(\E)$ is canonically unobstructed w.r.t. $\{\phi (t)\}_{t\in B}$.
\end{theorem}
\begin{proof}This follows from $4., 5., 6.$ of Proposition \ref{prop-12345}. \end{proof}
\section{Representing Dolbeault cohomology by invariant tensor fields} \label{Representing Dolbeault cohomology by invariant tensor fields}
By a \emph{holomorphic tensor bundle} on a complex manifold, we mean a holomorphic vector bundle formed by the tensor products or exterior products from the tangent bundle and the cotangent bundle. Let $M:=(\Gamma\backslash G,J)$ be a compact complex manifold where $G$ is a real Lie group with a lattice $\Gamma\subset G$ and $J$ is a left invariant complex structure on $G$. Given any holomorphic tensor bundle $E$ on $M$, there is a corresponding $\g^{0,1}$-module formed by replacing each tangent bundle component in $E$ by $\g^{1,0}$ and cotangent bundle component by $\g^{0,1}$. This $\g^{0,1}$-module will be denoted by $\E$.

Let $\pi: (\mathcal{M}, M)\to (B,0)$ be a small deformation of $M$ such that each fiber $M_t$ of $\pi$ is represented by an $G$-invariant $\phi (t)\in \Lambda^{0,1}\ot\g^{1,0}$ and $E_t$ is the holomorphic tensor bundle on $M_t$ corresponding to $E$. For each $t\in B$ the complex $(\Lambda^{0,\bullet}\ot \E, \pb_{\phi(t)})$ consisting of $G$-invariant tensor fields is naturally a subcomplex of
$(A^{0,\bullet}(M, E), \pb_{\phi(t)})$ and choose a $G$-invariant Hermitian metric on $M$, we have an orthogonal direct sum decomposition
\[
A^{0,\bullet}(M, E)=(\Lambda^{0,\bullet}\ot \E)\oplus (\Lambda^{0,\bullet}\ot \E)^\perp,
\]
such that $\left((\Lambda^{0,\bullet}\ot \E)^\perp, \pb_{\phi(t)}\right)$ is also a subcomplex. The $\pb_{\phi(t)}$-Laplacian operator $\Box_{\phi(t)}$ is defined as
\[
\Box_{\phi(t)}:=\pb_{\phi(t)}\pb_{\phi(t)}^*+\pb_{\phi(t)}^*\pb_{\phi(t)},
\]
and we have
\[
A^{0,\bullet}(M, E)\cap\ker\pb_{\phi(t)} =\mathcal{H}_{\pb_{\phi(t)}}^{0,\bullet}(M, E)\oplus (A^{0,\bullet}(M, E)\cap\Image\pb_{\phi(t)}),
\]
where $\mathcal{H}_{\pb_{\phi(t)}}^{0,\bullet}(M, E):=\ker\Box_{\phi(t)}\cap A^{0,\bullet}(M, E)$. It follows that
\[
\mathcal{H}_{\pb_{\phi(t)}}^{0,\bullet}(M, E)=\mathcal{H}^{0,\bullet}_{\pb_{\phi(t)}}(\E)\oplus \mathcal{H}^{0,\bullet}_{\pb_{\phi(t)}}(\E)^\perp,
\]
where $\mathcal{H}^{0,\bullet}_{\pb_{\phi(t)}}(\E):=\ker\Box_{\phi(t)}\cap (\Lambda^{0,\bullet}\ot \E)$ and $\mathcal{H}^{0,\bullet}_{\pb_{\phi(t)}}(\E)^\perp:=\ker\Box_{\phi(t)}\cap (\Lambda^{0,\bullet}\ot \E)^\perp$. It is easy to see that
\[
\ker\bar{\partial}_{\phi(t)} \cap (\Lambda^{0,\bullet}\ot \E)^\perp
=\mathcal{H}^{0,\bullet}_{\pb_{\phi(t)}}(\E)^\perp\oplus \left(\Image\bar{\partial}_{\phi(t)} \cap (\Lambda^{0,\bullet}\ot \E)^\perp\right) ,
\]
and
\[
\frac{\ker\bar{\partial}^*\cap\ker\bar{\partial}_{\phi(t)} \cap (\Lambda^{0,\bullet}\ot \E)^\perp}{\ker\bar{\partial}^*\cap\Image\bar{\partial}_{\phi(t)} \cap (\Lambda^{0,\bullet}\ot \E)^\perp}\cong \mathcal{H}^{0,\bullet}_{\pb_{\phi(t)}}(\E)^\perp~.
\]
\begin{theorem}\label{thm-main-result}
Let $M:=(\Gamma\backslash G,J)$ be a compact complex manifold where $G$ is a real Lie group with a lattice $\Gamma\subset G$ and $J$ is a left invariant complex structure on $G$. Denote by $\g$ the Lie algebra of $G$. Assume $E$ is a holomorphic tensor bundle on $M$. Let $\pi: (\mathcal{M}, M)\to (B,0)$ be a small deformation of $M$ such that each fiber $M_t$ of $\pi$ is represented by a $G$-invariant $\phi (t)\in \Lambda^{0,1}\ot\g^{1,0}$ and $E_t$ is the holomorphic tensor bundle on $M_t$ corresponding to $E$, then the set
\[
\left\{t\in B\mid H^{0,\bullet}_{\pb_t}(\E_t)\cong H^{0,\bullet}_{\pb_t}(M_t, E_t)\right\}
\]
is an analytic open subset (i.e. complement of analytic subset) of $B$.
\end{theorem}
\begin{proof}
It follows from the two extension isomorphisms (for the later one, see \cite[Thm.\,4.4]{Xia19dDol})
\[
H^{0,\bullet}_{\pb_t}(\E_t)\cong H^{0,\bullet}_{\pb_{\phi(t)}}(\E),\quad H^{0,\bullet}_{\pb_t}(M_t, E_t)\cong H^{0,\bullet}_{\pb_{\phi(t)}}(M, E)
\]
and the fact $\dim\mathcal{H}^{0,\bullet}_{\pb_{\phi(t)}}(\E)^\perp=\dim H^{0,\bullet}_{\pb_{\phi(t)}}(M, E)-\dim H^{0,\bullet}_{\pb_{\phi(t)}}(\E)$ that
\begin{align*}
\left\{t\in B\mid H^{0,\bullet}_{\pb_t}(\E_t)\cong H^{0,\bullet}_{\pb_t}(M_t, E_t)\right\}
=&\left\{t\in B\mid \dim\mathcal{H}^{0,\bullet}_{\pb_{\phi(t)}}(\E)^\perp=0\right\}\\
=&B\setminus\left\{t\in B\mid \dim\mathcal{H}^{0,\bullet}_{\pb_{\phi(t)}}(\E)^\perp\geq 1\right\}
\end{align*}
Hence it is enough to show $\{t\in B\mid \dim\mathcal{H}^{0,\bullet}_{\pb_{\phi(t)}}(\E)^\perp\geq 1\}$ is an analytic subset of $B$ which can be proved in the same way as Theorem \ref{thm-g-analytic-subset}.

In fact, we have
\begin{align*}
&\ker\bar{\partial}^*\cap\ker\bar{\partial}_{\phi(t)} \cap (\Lambda^{0,q}\ot \E)^\perp \cong\tilde{V}_{t}^q\\
:=&\{ \sum_{l=1}^N a_l\sigma_0^l\in \mathcal{H}^{0,q}_{\pb}(\E)^\perp \mid  (a_1, \cdots, a_N)\in \mathbb{C}^{N}~\text{s.t.}~\bar{\partial}_{\phi(t)}\sigma(t)=0,\\
&\text{where}~\sigma (t)=\sum_{k} \sigma_k~\text{with}~ \sigma_0=\sum_l a_l\sigma_0^l~\text{and}~ \sigma_{k}=\bar{\partial}^*G\sum_{i+j=k} \langle\phi_i| \sigma_{j} \rangle,~\forall k\neq 0 \}\\
=&\{ \sum_{l=1}^N a_l\sigma_0^l\in \mathcal{H}^{0,q}_{\pb}(\E)^\perp \mid \sum_{l=1}^N a_l \mathcal{H}^{0,q+1}\langle\phi(t) | \sigma^l(t) \rangle =0\}~,
\end{align*}
where $\sigma^l(t)=\sum_{k} \sigma_k^l$ with $\sigma_{k}^l=\bar{\partial}^\dag\sum_{i+j=k} \langle\phi_i| \sigma_{j}^l \rangle$ for $k> 0$.
We see that $\{t\in B\mid \dim \tilde{V}_{t}^q\geq k\}$ is an analytic subset of $B$ for any nonnegative integer $k$.

On the other hand, it can be shown similarly (see \cite[Prop.\,6.5]{Xia19dDol} and \cite[Thm.\,1.2]{Xia19dDol}) that
\[
\dim\mathcal{H}^{0,q}_{\pb_{\phi(t)}}(\E)^\perp=\dim\tilde{V}_{t}^q+\dim\tilde{V}_{t}^{q-1}-\dim\mathcal{H}^{0,q-1}_{\pb}(\E)^\perp.
\]
As a result, let $\lambda:=\dim \mathcal{H}^{0,q-1}_{\pb}(\E)^\perp$ we have
\begin{align*}
&\{t\in B\mid \dim\mathcal{H}^{0,q}_{\pb_{\phi(t)}}(\E)^\perp\geq 1\}\\
=&\{t\in B\mid \dim\tilde{V}_{t}^q+\dim\tilde{V}_{t}^{q-1}\geq 1+\lambda \}\\
=&\bigcup_{j=0}^\lambda\Big(\{t\in B\mid \dim\tilde{V}_{t}^q\geq 1+\lambda-j\}\cap\{t\in B\mid \dim\tilde{V}_{t}^{q-1}= j\}\Big)\\
=&\bigcup_{j=0}^\lambda\Big(\{t\in B\mid \dim\tilde{V}_{t}^q\geq 1+\lambda-j\}\cap\{t\in B\mid \dim\tilde{V}_{t}^{q-1}\geq j\}\Big).
\end{align*}
The conclusion follows from the fact that arbitrary intersections or finite unions of analytic subsets are still analytic subsets (see e.g. ~\cite[pp.86-87]{GR65}).
\end{proof}

For the Bott-Chern cohomology on Lie algebra with complex structures, we have a similar theory as those developed in \cite{Xia19dBC}. The analogue of Theorem \ref{thm-main-result} still holds in this case. This justifies our computations of the dimensions of deformed Bott-Chern cohomology in \cite{Xia19dBC}.

Now we present some explicit computations.
\begin{example}
Let $X=\mathbb{C}^3/\Gamma$ be the solvable manifold constructed by Nakamura in Example III-(3b) of~\cite{Nak75}. We have
\begin{align*}
H^0(X,\Omega_X^{1})~ =& ~\mathbb{C}\{\varphi^1 = d z^1,~ \varphi^2 = e^{z_1} d z^2,~ \varphi^3 = e^{-z_1} d z^3\} ~,\\
H^0(X,T_X^{1,0})~  =& ~\mathbb{C}\{\theta^1=\frac{\partial}{\partial z^1},~\theta^2= e^{-z_1} \frac{\partial}{\partial z^2},~\theta^3= e^{z_1} \frac{\partial}{\partial z^3}\} ~,\\
\mathcal{H}^{0,1}(X)~ =& ~\mathbb{C}\{ \psi^{\bar{1}} = d z^{\bar{1}},~ \psi^{\bar{2}} = e^{z_1} d z^{\bar{2}},~ \psi^{\bar{3}} = e^{-z_1} d z^{\bar{3}}\} ~,\\
\mathcal{H}^{0,1}(X, T_X^{1,0})~ =& ~\mathbb{C}\{\theta^i\psi^{\bar{\lambda}} ,~ i=1, 2, 3, \lambda=1, 2, 3\} ~,
\end{align*}
where we use the Hermitian metric $\sum_{i=1}^3\varphi^i\otimes\bar{\varphi}^i$. The Beltrami differential for the fiber $X_t$ of the Kuranishi family of $X$ is given by
\[
\phi(t) = \phi_1 = t_{i\lambda}\theta^i\psi^{\bar{\lambda}}
\]
and the Kuranishi space of $X$ is
\[
\mathcal{B}=\{t=(t_{11}, t_{12}, t_{13}, t_{21}, t_{22}, t_{23}, t_{31}, t_{32}, t_{33})\in \mathbb{C}^9\mid |t_{i\lambda}|<\epsilon, i=1, 2, 3, \lambda=1, 2, 3 \},
\]
where $\epsilon>0$ is sufficiently small. By Theorem \ref{thm-main-result} and Theorem \ref{thm-ext-iso} respectively, we have for any $t\in\mathcal{B}$ there holds $\dim H^{p,q}_{\pb_t}(X_t)=\dim H^{p,q}_{\pb_t}(\g, J_t)=\dim H^{p,q}_{\pb_{\phi(t)}}(\g, J)$ which can be determined by analysing the deformation behavior of classes in $H^{p,q}_{\pb}(\g, J)$ (See Proposition \ref{prop-12345}).

Let us now consider the deformation of classes in $H^{1,0}_{\pb}(\g, J)$ and set $\sigma_0 = a_{1}\varphi^{1} +a_{2}\varphi^{2} +a_{3}\varphi^{3} $, then
\begin{align*}
\mathcal{L}_{\phi_1}^{1,0}\sigma_0
=& a_{1}(-t_{12}\varphi^{1}\wedge\psi^{\bar{2}} + t_{13}\varphi^{1}\wedge\psi^{\bar{3}} ) + a_{2}(-t_{1\lambda}\varphi^{2}\wedge\psi^{\bar{\lambda}} + t_{21}\varphi^{1}\wedge\psi^{\bar{1}} + 2t_{23}\varphi^{1}\wedge\psi^{\bar{3}})\\
 &+ a_{3}( t_{1\lambda}\varphi^{3}\wedge\psi^{\bar{\lambda}} - t_{31}\varphi^{1}\wedge\psi^{\bar{1}} - 2t_{32}\varphi^{1}\wedge\psi^{\bar{2}} )\\
=& (a_2t_{21}-a_3t_{31})\varphi^{1}\wedge\psi^{\bar{1}} - (a_1t_{12} + 2a_3t_{32})\varphi^{1}\wedge\psi^{\bar{2}} + (a_1t_{13} + 2a_2t_{23})\varphi^{1}\wedge\psi^{\bar{3}}\\
 &+ a_{3}t_{1\lambda}\varphi^{3}\wedge\psi^{\bar{\lambda}} - a_{2}t_{1\lambda}\varphi^{2}\wedge\psi^{\bar{\lambda}}
\end{align*}
is exact if and only if $\mathcal{L}_{\phi_1}^{1,0}\sigma_0=0$, i.e.
\begin{equation}\label{a_{2}}
\left\{
\begin{array}{rcl}
a_{2}t_{21}- a_{3}t_{31} &=& 0 \\[5pt]
a_{1}t_{12}+ 2a_{3}t_{32} &=& 0 \\[5pt]
a_{1}t_{13}+ 2a_{2}t_{23} &=& 0 \\[5pt]
a_{2}t_{1\lambda} &=& 0, ~ \lambda=1,~2,~3 \\[5pt]
a_{3}t_{1\lambda} &=& 0, ~ \lambda=1,~2,~3 ~.
\end{array}
\right.
\end{equation}
has solutions for $(a_{1}, a_{2}, a_{3})$. On the other hand, $\phi_k=0,~ k>1$  implies that $\sigma_{k}=0,~ k>1$.

Therefore, for $V=H^{1,0}_{\pb}(\g, J)$ we have
\begin{align*}
V_{t}^0=
&\{ a_{1}\varphi^{1} +a_{2}\varphi^{2} +a_{3}\varphi^{3} \mid  (a_{1}, a_{2}, a_{3})\in \mathbb{C}^3~\text{s.t.}~\bar{\partial}_{\phi(t)}\sigma(t)=0,\\
&\text{where}~\sigma (t)=\sum_{k} \sigma_k~\text{with}~ \sigma_0=\sum_l a_l\sigma_0^l~\text{and}~ \sigma_{k}=\bar{\partial}^\dag\sum_{i+j=k} \langle\phi_i| \sigma_{j} \rangle,~\forall k\neq 0 \}\\
=&\{a_{1}\varphi^{1} +a_{2}\varphi^{2} +a_{3}\varphi^{3} \mid (a_{1}, a_{2}, a_{3})\in \mathbb{C}^3~\text{satisfy}~ \eqref{a_{2}}\}
\end{align*}
and $\dim V_t^0$ is determined by the rank of the coefficient matrix $T$ of \eqref{a_{2}}.

Hence we have the following table
\vspace{12pt}
\renewcommand\arraystretch{1.5}
\begin{table}[!htbp]
\centering
\begin{center}
\begin{tabular}{|c|c|c|c|}
\hline
$t\in \mathcal{B}$ & $\text{rank}~ T$ & $\dim V_t^0$ & $\dim H^{1,0}_{\pb_t}(X_t)$ \\
\hline
$t_{i\lambda}=0, 1\leq i,\lambda\leq 3, (i,\lambda)\neq (2,2), (3,3)$ & $0$ & $3$ & $3$   \\
\hline
\tabincell{c}{$t_{11}=t_{12}=t_{13}=t_{21}t_{32}=t_{23}t_{31}=t_{23}t_{32}=0$ s.t.\\ $t_{23}\neq 0$ or $t_{32}\neq 0$ or $t_{21}\neq 0$ or $t_{31}\neq 0$} & $1$ & $2$ & $2$ \\
\hline
\tabincell{c}{$t_{12}=t_{13}=0$ s.t.\\ $t_{11}\neq 0$ or $t_{23}t_{32}\neq 0$ or $t_{23}t_{31}\neq 0$ or $t_{21}t_{32}\neq 0$} & $2$ & $1$ & $1$ \\
\hline
$t_{12}\neq 0$ or $t_{13}\neq 0$ or $t_{21}t_{13}t_{32}-t_{31}t_{12}t_{23}\neq 0$ & $3$ & $0$ & $0$ \\
\hline
\end{tabular}
\end{center}
\end{table}
\end{example}

\vskip 1\baselineskip \textbf{Acknowledgements.} I would like to thank Prof. Kefeng Liu for constant encouragement and Prof. Nailin Du for introducing me to the theory of Moore-Penrose inverses. Many thanks to S\"onke Rollenske, Mutaz Abumathkur and Jialin Zhu for useful communications.

\bibliographystyle{alpha}

\begin{thebibliography}{CFGU00}

\bibitem[AK17]{AK17b}
D.~Angella and H.~Kasuya.
\newblock Cohomologies of deformations of solvmanifolds and closedness of some
  properties.
\newblock {\em North-West. Eur. J. Math.}, 3:75--105, 2017.

\bibitem[Ang13]{Ang13}
D.~Angella.
\newblock The cohomologies of the {Iwasawa} manifold and of its small
  deformations.
\newblock {\em J. Geom. Anal.}, 23(3):1355--1378, 2013.

\bibitem[BDV09]{BDV09}
M.~L. Barberis, I.~G. Dotti, and M.~Verbitsky.
\newblock Canonical bundles of complex nilmanifolds, with applications to
  hypercomplex geometry.
\newblock {\em Math. Res. Lett.}, 16(2):331--347, 2009.

\bibitem[CF01]{CF01}
S.~Console and A.~Fino.
\newblock Dolbeault cohomology of compact nilmanifolds.
\newblock {\em Transform. Groups}, 6(2):111--124, 2001.

\bibitem[CFGU00]{CFGU00}
L.~A. Cordero, M.~Fern\'{a}ndez, A.~Gray, and L.~Ugarte.
\newblock Compact nilmanifolds with nilpotent complex structures: {D}olbeault
  cohomology.
\newblock {\em Trans. Amer. Math. Soc.}, 352(12):5405--5433, 2000.

\bibitem[CFK16]{CFK16}
S.~Console, A.~Fino, and H.~Kasuya.
\newblock On de {R}ham and {D}olbeault cohomology of solvmanifolds.
\newblock {\em Transform. Groups}, 21(3):653--680, 2016.

\bibitem[CFP06]{CFP06}
S.~Console, A.~Fino, and Y.~S. Poon.
\newblock Stability of abelian complex structures.
\newblock {\em Internat. J. Math.}, 17(4):401--416, 2006.

\bibitem[Con06]{Con06}
S.~Console.
\newblock Dolbeault cohomology and deformations of nilmanifolds.
\newblock {\em Rev. Un. Mat. Argentina}, 47(1):51--60, 2006.

\bibitem[FRR19]{FRR19}
A.~Fino, S.~Rollenske, and J.~Ruppenthal.
\newblock Dolbeault cohomology of complex nilmanifolds foliated in toroidal
  groups.
\newblock {\em Q. J. Math.}, 70(4):1265--1279, 2019.

\bibitem[GR65]{GR65}
R.~C. Gunning and H.~Rossi.
\newblock {\em Analytic functions of several complex variables}.
\newblock Prentice-Hall, Inc., Englewood Cliffs, N.J., 1965.

\bibitem[Gro77]{Gro77}
C.~W. Groetsch.
\newblock {\em Generalized inverses of linear operators: representation and
  approximation}.
\newblock Marcel Dekker, Inc., New York-Basel, 1977.
\newblock Monographs and Textbooks in Pure and Applied Mathematics, No. 37.

\bibitem[GT93]{GT93}
G.~Gigante and G.~Tomassini.
\newblock Deformations of complex structures on a real {L}ie algebra.
\newblock In {\em Complex analysis and geometry}, Univ. Ser. Math., pages
  377--385. Plenum, New York, 1993.

\bibitem[Kas13]{Kas13}
H.~Kasuya.
\newblock Techniques of computations of {D}olbeault cohomology of
  solvmanifolds.
\newblock {\em Math. Z.}, 273(1-2):437--447, 2013.

\bibitem[Kas16]{Kas16}
H.~Kasuya.
\newblock An extention of {N}omizu's {T}heorem--a user's guide.
\newblock {\em Complex Manifolds}, 3(1):231--238, 2016.

\bibitem[KMS93]{KMS93}
I.~Kol{\'a}$\check{r}$, P.~W. Michor, and J.~Slov{\'a}k.
\newblock {\em Natural operations in differential geometry}.
\newblock Springer-Verlag, Berlin, Heidelberg, New York, 1993.

\bibitem[LR11]{LR11}
K.~Liu and S.~Rao.
\newblock Remarks on the {C}artan formula and its applications.
\newblock {\em Asian J. Math.}, 16(1):157--169, 2011.

\bibitem[LRY15]{LRY15}
K.~Liu, S.~Rao, and X.~Yang.
\newblock Quasi-isometry and deformations of {C}alabi-{Y}au manifolds.
\newblock {\em Invent. Math.}, 199(2):423--453, 2015.

\bibitem[LSY09]{LSY09}
K.~Liu, X.~Sun, and S.-T. Yau.
\newblock Recent development on the geometry of the {Teichm\"uller} and moduli
  spaces of {Riemann} surfaces.
\newblock In {\em Geometry of Riemann surfaces and their moduli spaces}, volume
  XIV of {\em {Surveys in differential geometry}}, pages 221--259. 2009.

\bibitem[LZ18]{LZ18}
K.~Liu and S.~Zhu.
\newblock Solving equations with {Hodge} theory.
\newblock arXiv:1803.01272v1, 2018.

\bibitem[MK06]{MK71}
J.~Morrow and K.~Kodaira.
\newblock {\em Complex manifolds}.
\newblock AMS Chelsea Publishing, Providence, RI, 2006.
\newblock Reprint of the 1971 edition with errata.

\bibitem[MPPS06]{MPPS06}
C.~Maclaughlin, H.~Pedersen, Y.~S. Poon, and S.~Salamon.
\newblock Deformation of 2-step nilmanifolds with abelian complex structures.
\newblock {\em J. London Math. Soc. (2)}, 73(1):173--193, 2006.

\bibitem[Nak75]{Nak75}
I.~Nakamura.
\newblock Complex parallelisable manifolds and their small deformations.
\newblock {\em J. Differential Geom.}, 10(1):85--112, 1975.

\bibitem[Nom54]{Nom54}
K.~Nomizu.
\newblock On the cohomology of compact homogeneous spaces of nilpotent {L}ie
  groups.
\newblock {\em Ann. of Math. (2)}, 59:531--538, 1954.

\bibitem[OV20]{OV20}
L.~Ornea and M.~Verbitsky.
\newblock Twisted {D}olbeault cohomology of nilpotent lie algebras.
\newblock {\em Transform. Groups}, 2020.
\newblock https://doi.org/10.1007/s00031-020-09601-4.

\bibitem[Rol09a]{Rol09large}
S.~Rollenske.
\newblock Geometry of nilmanifolds with left-invariant complex structure and
  deformations in the large.
\newblock {\em Proc. Lond. Math. Soc. (3)}, 99(2):425--460, 2009.

\bibitem[Rol09b]{Rol09}
S.~Rollenske.
\newblock Lie-algebra {D}olbeault-cohomology and small deformations of
  nilmanifolds.
\newblock {\em J. Lond. Math. Soc. (2)}, 79(2):346--362, 2009.

\bibitem[Rol11]{Rol11b}
S.~Rollenske.
\newblock Dolbeault cohomology of nilmanifolds with left-invariant complex
  structure.
\newblock In {\em Complex and differential geometry}, volume~8 of {\em Springer
  Proc. Math.}, pages 369--392. Springer, Heidelberg, 2011.

\bibitem[RTW20]{RTW20}
S.~Rollenske, A.~Tomassini, and X.~Wang.
\newblock Vertical-horizontal decomposition of {L}aplacians and cohomologies of
  manifolds with trivial tangent bundles.
\newblock {\em Ann. Mat. Pura Appl. (4)}, 199(3), 2020.

\bibitem[RWZ19]{RWZ19}
S.~Rao, X.~Wan, and Q.~Zhao.
\newblock On local stabilities of $p$-{K\"ahler} structures.
\newblock {\em Compos. Math.}, 155(3):455--483, 2019.

\bibitem[RZ18]{RZ18}
S.~Rao and Q.~Zhao.
\newblock Several special complex structures and their deformation properties.
\newblock {\em J. Geom. Anal.}, 28(4):2984--3047, 2018.

\bibitem[Wei94]{Wei94}
C.~Weibel.
\newblock {\em An introduction to homological algebra}, volume~38 of {\em
  Cambridge Studies in Advanced Mathematics}.
\newblock Cambridge University Press, Cambridge, 1994.

\bibitem[WWQ18]{WWQ18}
G.~Wang, Y.~Wei, and S.~Qiao.
\newblock {\em Generalized inverses: theory and computations}, volume~53 of
  {\em Developments in Mathematics}.
\newblock Springer, Singapore; Science Press Beijing, Beijing, second edition,
  2018.

\bibitem[Xia19a]{Xia19dDol}
W.~Xia.
\newblock Deformations of {Dolbeault} cohomology classes.
\newblock arXiv:1909.03592, 2019.

\bibitem[Xia19b]{Xia19deri}
W.~Xia.
\newblock Derivations on almost complex manifolds.
\newblock {\em Proc. Amer. Math. Soc.}, 147:559--566, 2019.
\newblock Errata in arXiv:1809.07443v3.

\bibitem[Xia20]{Xia19dBC}
W.~Xia.
\newblock On the deformed {Bott-Chern} cohomology.
\newblock {\em J. Geom. Phys.}, 166:104250, 2021.

\bibitem[ZR13]{ZR13}
Q.~Zhao and S.~Rao.
\newblock Applications of the deformation formula of holomorphic one-forms.
\newblock {\em Pacific J. Math.}, 266(1):221--255, 2013.

\bibitem[ZR15]{RZ15}
Q.~Zhao and S.~Rao.
\newblock Extension formulas and deformation invariance of {Hodge} numbers.
\newblock {\em C. R. Math. Acad. Sci. Paris}, 353(11):979--984, 2015.

\end{thebibliography}

\end{document}